
\documentstyle{amsppt}
\pagewidth{6.4in}\vsize8.5in\parindent=6mm\parskip=3pt\baselineskip=14pt\tolerance=10000 \hbadness=500
\loadbold
\topmatter
\title
Bounds for singular fractional integrals  and
 related Fourier integral  operators
\endtitle
\author
Andreas Seeger \ \ \ \ \ \ \ \ \ \ \ \ \ \ \ \ \  Stephen Wainger
\endauthor
\rightheadtext{Bounds for singular fractional integrals}
\thanks 
Research supported in part by  NSF grants DMJ 9731647 and DMJ 9970042.
\endthanks
\date September, 15 2000\enddate
\address 
Department of Mathematics, University of Wisconsin-Madison, Madison, WI 53706
\endaddress 
\email seeger\@math.wisc.edu\endemail
\address 
Department of Mathematics, University of Wisconsin-Madison, Madison, WI 53706
\endaddress 
\email wainger\@math.wisc.edu\endemail
\endtopmatter
\document


\define\rank{{\text{\rm rank }}}

\define\supp{{\text{\rm supp }}}

\define\inn#1#2{\langle#1,#2\rangle}

\define\noi{\noindent}
\define\lcontr{\rfloor}
\define\lco#1#2{{#1}\lcontr{#2}}
\define\lcoi#1#2{\imath({#1}){#2}}
\define\rco#1#2{{#1}\rcontr{#2}}

\define\bin#1#2{{\pmatrix {#1}\\{#2}\endpmatrix}}

\define\card{\text{\rm card}}
\define\lc{\lesssim}
\define\gc{\gtrsim}


\define\eps{\varepsilon}

\define\ka{\kappa}
            
\define\la{\lambda}

\define\om{\omega}              \define\Om{\Omega}

\define\fM{{\frak M}}

\define\fQ{{\frak Q}}

\define\fS{{\frak S}}


\define\bbR{{\Bbb R}}
\define\bbS{{\Bbb S}}

\define\bbZ{{\Bbb Z}}

\define\cA{{\Cal A}}
\define\cB{{\Cal B}}
\define\cC{{\Cal C}}

\define\cE{{\Cal E}}

\define\cG{{\Cal G}}

\define\cI{{\Cal I}}
\define\cJ{{\Cal J}}

\define\cL{{\Cal L}}
\define\cM{{\Cal M}}

\define\cQ{{\Cal Q}}
\define\cR{{\Cal R}}
\define\cS{{\Cal S}}
\define\cT{{\Cal T}}

\define\cV{{\Cal V}}




\define\r{\rho}
\define\s{\sigma}
\define\vth{\vartheta}
\define\vect{{\text{vect}}}

\define\ic{\imath}
\define\innmod#1#2{#1\!\cdot\!#2}        
\redefine\Re{\text{Re}}

\head{\bf 1. Introduction}\endhead

Let $\Omega\subset \widetilde \Omega $  be   open sets in $\Bbb R^d$, $I\subset \bbR^{d-\ell}$  
be an open neighborhood of
the origin and  let $ \eta$ be a compactly 
supported  smooth function  on $\Omega\times I$; we
assume that
$\eta(\cdot,0)$ does not vanish identically. For each $x\in \Omega$ let
$t\mapsto \Gamma(x,t) \subset\widetilde \Omega $ be a regular parametrization 
of a 
submanifold $\cM_x\subset\widetilde \Omega$ with codimension $\ell$.
We assume that $\Gamma(x,t)\subset \Omega$ if $(x,t)\in\supp \eta$, and that
$\Gamma$ satisfies
$\Gamma(x,0)=x$ and  depends smoothly on $(x,t)$.

We shall consider the singular  fractional integral operator
(or  weakly singular Radon transform)
$\cR_\sigma$, defined by
$$\cR_\sigma f(x)
= \int \eta(x,t) f(\Gamma(x,t)) 
|t|^{-(d-\ell-\sigma)} \, dt,
\tag 1.1
$$
under suitable ``curvature'' assumptions on the singular support 
and the wavefront sets of the distribution kernel of the integral operator.

To formulate these assumptions  we shall 
work with a submanifold $\cM$ 
of codimension $\ell$ in $ \Omega\times\Omega$; 
so that 
$$\Delta:=\{(x,x): x\in \Omega\}\subset \cM.
\tag 1.2
$$
To relate this to the operator in (1.1) we assume that $\eta(x,t)$ 
vanishes unless $|t|<\delta$ for small $\delta$ and note that 
 the differential of the map $(x,t)\mapsto \gamma(x,t)$ 
has maximal rank $d+\ell$; then we take
$$\cM=\{(x,y): x\in \Omega, y=\Gamma(x,t) \text{ for  some } |t|<\delta\}.
$$

Moreover we assume the following standard hypotheses in the theory of Fourier
integral operators:

\proclaim{ { Nondegeneracy assumptions}}

\noi {$(1.3)$} The natural projections
$(x,y)\mapsto x$ and $(x,y)\mapsto y$ are submersions when restricted to 
$\cM$.

\noi {$(1.4)$} The twisted normal bundle
$N^*\cM'\subset T^*\Om\times T^*\Om$
is locally the graph of a canonical transformation. 
Here $N^*\cM'$ consists of all $(x,\xi,y,-\eta)$ where
$(x,y)\in \cM$ and $(\xi,\eta)\in T_{(x,y)}^*\cM$ annihilates the tangent vectors
in $T_{(x,y)}\cM$.
\endproclaim

Assumption (1.3) implies  that the sections
$$
\align
\cM_x&=\{y\in \Omega: (x,y)\in \cM\}
\\
\cM^y&=\{x\in \Omega: (x,y)\in \cM\}
\endalign
$$
are immersed submanifolds of $\Om$, of codimension $\ell$.
We may assume that 
 $\cM$ is given by a defining  function
$$\cM=\{(x,y):\Phi(x,y)=0\}
\tag 1.5$$ 
where
$\Phi$ is  $\Bbb R^\ell$-valued satisfying $\Phi(x,x)=0$ so that (1.2) is satisfied
and $\rank \Phi_x=\rank\Phi_y=\ell$ so that (1.3) is satisfied.

Assumption (1.4) can be reformulated  as follows. 
Let  $\Psi(x,y,\tau)=\tau\!\cdot\!\Phi(x,y)$. Then the 
assumption (1.4) on $N^*\cM'$ 
is  equivalent with
$$
\det \pmatrix
\Psi_{xy}&\Psi_{x\tau}
\\
\Psi_{\tau y}&\Psi_{\tau\tau}\endpmatrix
=
\det\pmatrix
\tau\!\cdot\!\Phi_{xy}&\Phi_x
\\
\Phi_y&0\endpmatrix \neq 0
\quad\text{ for all $\tau \in S^{\ell-1}$, }
\tag 1.6
$$
see \cite{12}; in (1.6) 
$\Phi_x$ should be read as a $d\times\ell$ - matrix and $\Phi_y $ as an $\ell\times d$ -  matrix.
For   $\ell=1$ hypothesis (1.4) is just the  rotational 
curvature assumption of Phong and Stein \cite{16}.
We note that for (1.4) to hold the codimension $\ell$ 
has to be sufficiently small, 
and we are mainly interested in the case of hypersurfaces.


\proclaim{Theorem 1.1} Suppose that $1\le p\le  q\le \infty$,  $0<\s<d-\ell$ 
and suppose that $\cM$  satisfies the nondegeneracy assumptions (1.3), (1.4).
Then 
$\cR_\sigma$ maps $L^p(\Omega)\to L^q (\widetilde \Omega)$
if and only if  the following conditions are satisfied:

\roster
\item"{{\it a})}"
 $(1/p, 1/q)$ belongs to the triangle with corners
$(0,0)$, $(1,1)$ and $(\frac{d}{d+\ell}, \frac{\ell}{d+\ell})$.

\item"{{\it b})}"
 $(1/p, 1/q)$ belongs to the halfplane defined by
 $(d+\ell)(\frac 1p-\frac 1q)\le \s$.
\endroster 
\endproclaim

A special translation invariant case is due to M. Christ \cite{3},
extending earlier results by Ricci and Stein \cite{18}. These authors consider
 the translation invariant  case where 
 $\Phi(x,y)=x_d-y_d-|x'-y'|^2$ and 
a related model case on the Heisenberg group. For these dilation invariant
examples 
 one actually 
 proves global results which one could deduce from 
 local  ones by scaling arguments.

The weakly singular Radon transforms are special cases of 
oscillatory integrals  with singular symbols as considered by Melrose
\cite{13}, Greenleaf and Uhlmann \cite{11} and others. Let 
$I^{\r,-\s}(\Omega\times\Omega;\cM,\Delta)$ denote 
the class of distribution kernels introduced in \cite{11}; we denote by
$\cI^{\r,-\s}(\Omega\times\Omega;\cM,\Delta)$
 the associated class of operators
and refer for a general discussion and other references to previous work 
to \cite{11}.

Possibly after a change of variable we may locally  parametrize 
 $\cM$  as a graph  of an $\Bbb R^\ell$ valued function,
$$
y''=S(x,y')
\tag 1.7
$$
with $y'=(y_1,\dots, y_{d-\ell})$, $y''=(y_{d-\ell+1},\dots, y_d)$,
$S=(S_{d-\ell+1},\dots, S_d)$. 
so that
$$\rank S_{x''} =\ell\tag 1.8$$
and
$$\det \pmatrix
\innmod{\theta}{S_{x'y'}} &S_{x'}
\\
\innmod{\theta}{S_{x''y'}} &S_{x''}
\endpmatrix \neq 0
\tag 1.9
$$
for all $\theta\in \bbR^\ell\setminus\{0\}$.

We recall from \cite{11}, \cite{5} that a  distribution kernel $K$ belongs to
$ I^{\r,-\s}(\Omega\times\Omega;\cM,\Delta)$ if it is  
a locally finite  sum of  $K_\nu$, so that each $K_\nu$
can be written after a change of variable in $\Omega$ as an oscillatory 
integral
$$
\iint_{\Bbb R^{d-\ell}\times \Bbb R^\ell}
 e^{\imath[\inn{\tau}{y''-S(x,y')}+\inn{\xi}{x'-y'}]}
a(x,y,\tau,\xi) d\tau d\xi.
\tag 1.10
$$
Here $S$ satisfies (1.8), (1.9) and the symbol 
 $a$ satisfies the differential inequalities
$$
|\partial_{x,y}^\gamma \partial_\tau^\alpha\partial_\xi^\beta
a(x,y,\tau,\xi)|\le C_{\alpha, \beta,\gamma}
(1+|\tau|+|\xi|)^{\r-|\alpha|} (1+|\xi|)^{-\s-|\beta|}.
\tag 1.11
$$
We refer to the class of symbols satisfying (1.11) as
$S^{\r,-\s}(\Omega\times\Omega,\bbR^{\ell}, \bbR^{d-\ell})$.
We shall sometimes  denote the operator with kernel (1.10) as $\cT[a]$.

It is well known  that the weakly singular Radon transform 
as considered in Theorem 1.1 is  an
 operator in 
$ \cI^{0,-\s}(\Omega\times\Omega;\cM,\Delta)$ (see {\it e.g.} \cite{11}).
Namely 
after an  appropriate localization it suffices to work with
$$
\cR^\s f(x)=\int f(y', S(x,y')) |x'-y'|^{\sigma+\ell-d}\chi(x', S(x,y'),y') dy'
\tag 1.12 
$$
where $\chi$ has small support. Then the distribution  kernel is given by
$$ \delta(y''-S(x,y'))|x'-y'|^{\s-d+\ell}g(x,y')$$
where  $\delta$ is the Dirac measure at the origin in  $\Bbb R^\ell$ and $g$ is
smooth and compactly supported.
We expand the Dirac measure using the   Fourier inversion formula in $\Bbb R^\ell$
 and  apply the Fourier inversion formula in $\Bbb R^{d-\ell}$
to the  function $h\to |h|^{\s-d+\ell}g(x,x'+h)$.
As a result we can write the distribution kernel  in the form
(1.10) where the symbol $a$  is given by
$$a(x,y,\tau,\xi)= (2\pi)^{-d} \int |w'|^{\s-d+\ell}
g(x,x'+w') e^{-\imath\inn {\xi} {w'}}dw'.
$$

We now formulate estimates for general operators of class
$\cI^{\r,-\s}$. 
Since the composition of a standard pseudo-differential operator
of order $m$ with an operator in 
$\cI^{\r, -\s}(\Omega\times\Omega;\cM,\Delta)$ belongs to
$\cI^{\r+m, -\s}(\Omega\times\Omega;\cM,\Delta)$  (see \cite{5}, \cite{11})
the following results  yield 
$L^p_a\to L^q_{a+m}$ Sobolev estimates for weakly singular Radon transforms.

\proclaim{Theorem 1.2} Suppose that $1\le p\le q\le\infty$.
Let $T\in \cI^{\r, -\s}(\Omega\times\Omega;\cM,\Delta)$,
with compactly supported distribution kernel, and assume that the nondegeneracy assumptions
(1.3), (1.4) hold.

\noi {\bf 1.2.1.} Suppose $0< \r <\frac{d-\ell}2$ and $2\r< \s<d-\ell$.
Then
$T$ maps $L^p$ to $L^q$ if 
the following two conditions are
 satisfied.

\roster
\item"{\it a)}"
 $(1/p, 1/q)$ belongs to the closed triangle with corners $(\frac{\r}{d-\ell}, \frac{\r}{d-\ell})$,
$(\frac{d-\ell-\r}{d-\ell},\frac{d-\ell-\r}{d-\ell})$ and 
$(\frac{d-\r}{d+\ell}, \frac{\r+\ell}{d+\ell})$.

\item"{\it b)}"
 $(1/p, 1/q)$ belongs to the halfspace defined by
 $(d+\ell)(\frac 1p-\frac 1q)\le 
\s-2\r$.
\endroster

\noi {\bf 1.2.2.} Suppose $\r=0$ and $0< \s<d-\ell$.
Then
$T$ maps $L^p$ to $L^q$ if the following two conditions are
 satisfied.

\roster
\item"{\it a)}"
 $(1/p, 1/q)$ belongs to the closed triangle with corners $(0,0)$, $(1,1)$ and 
$(\frac{d}{d+\ell}, \frac{\ell}{d+\ell})$, with the possible exception of the points
$(0,0)$ and $(1,1)$.

\item"{\it b)}"
 $(1/p, 1/q)$ belongs to the halfspace defined by
 $(d+\ell)(\frac 1p-\frac 1q)\le 
\s$.
\endroster

Moreover,  $T$ is bounded from the Hardy space $H^1$ to $L^1$ and from $L^\infty$ to $BMO$.

\noi {\bf 1.2.3.}
 Suppose $-\ell< \r< 0$ and  $-\r\frac{d-\ell}{\ell}< \s<d-\ell $.
Then
$T$ maps $L^p$ to $L^q$ if the following two conditions are
 satisfied.

\roster
\item"{\it a)}"
 $(1/p, 1/q)$ belongs to the pentagon with  corners
$(1,1)$, $(0,0)$, $(1,\frac{\rho+\ell}{\ell})$, $(\frac{-\rho}{\ell},0)$ and
$(\frac{d-\r}{d+\ell}, \frac{\r+\ell}{d+\ell})$, with the possible exceptions of the points
$(1,\frac{\rho+\ell}{\ell})$, $(\frac{-\rho}{\ell},0)$.

\item"{\it b)}"  $(1/p, 1/q)$ belongs to the halfspace defined by
 $(d+\ell)(\frac 1p-\frac 1q)\le 
\s-2\r$.
\endroster

\noi {\bf 1.2.4.}
 Suppose $-\ell< \r< 0$ and  $0<\s\le -\r\frac{d-\ell}{\ell} $. Then $T$ 
 maps $L^p$ to $L^q$ if
 $(1/p, 1/q)$ belongs to the quadrilateral with corners
$(1,1)$, $(0,0)$, $(1, \frac{d-\s-\ell}{d-\ell})$ and $(\frac{\s}{d-\ell},0)$,
with the possible exception of the points
$(1, \frac{d-\s-\ell}{d-\ell})$ and $(\frac{\s}{d-\ell},0)$.
\endproclaim

We remark that the analytic family of fractional integrals
 considered by Grafakos \cite{9} in the translation invariant case
can be considered  as a model family of operators of class $\cI^{\r,-\s}$, 
however 
the $L^2$ endpoint case in this family belongs to
$\cI^{\frac{d-\ell}2,\ell-d}$ but  satisfies better 
$L^2$ estimates than the general 
operator in $\cI^{\frac{d-\ell}2,\ell-d}$.

Operators in $\cI^{0,0}$ are bounded on $L^p$  for $1<p<\infty$,
see Greenleaf and Uhlmann \cite{11}, and for the main  special case of 
singular Radon transforms 
Phong and Stein \cite{16}, \cite{17}. 
The endpoint $L^p\to L^p$ estimates for the case $2\rho=\s$, 
$p_\rho=(d-\ell-\rho)/(d-\ell)$ or
$p_\rho'=\rho/(d-\ell)$ may fail  as demonstrated by Christ \cite{4}. 
It is  likely that 
the best possible  Lorentz-space endpoint estimate, namely an $L^{p_\rho}\to 
L^{p_\rho,2}$ bound holds;
a proof of this estimate  in the translation-invariant case was given by Tao and 
one of the authors \cite{20}.

A variant of the  methods in this paper has been used by the authors \cite{21}   to prove  new  $L^p$ theorems 
 for variable-coefficient 
maximal  and singular integral operators associated to families of curves 
in $\bbR^2$ (extending  results in \cite{2}, \cite{19}).

It is well known that at least under the assumption of
nonvanishing rotational curvature  
certain parabolic cutoffs can be used to write a singular integral
 along a hypersurface  as a sum of two
 operators, where  one of them
is  a  pseudodifferential operator of type $(1/2,1/2)$ and the
other one a Fourier integral operator, of type  $(1/2,1/2)$.
This decomposition  is due to Melrose (see \cite{13}, \cite{11}), but 
related arguments 
had been used by 
Nagel, Stein and Wainger \cite{15}, see also Phong and Stein \cite{17} for 
a different version.
In the course of this paper we shall make use of 
(variants of) all these decompositions.

The paper is organized as follows: 
\S2 contains some preparations and the discussion of a
 crucial change 
of variables. \S3 contains preliminary estimates for dyadic pieces of fractional Radon 
transforms. After appropriate localizations these are reduced to standard estimates for 
Fourier integral operators via
parabolic scalings.
In  $\S4$ we consider  some variants of  fractional integrals 
which are relevant for the estimation of the 
pseudodifferential contribution to 
 operators in $\cI^{\r,-\s}$ when $\rho\le 0$. Here we shall  also see that part 1.2.4 follows in a straightforward way from estimates for a class for certain product-type fractional integrals.
In \S5    we give the proof of Theorem 1.1. It turns out that after some  changes of 
variables  angular Littlewood-Paley decompositions may be applied just as in the previously known
translation-invariant case (\cite{3}). As in that case a positivity argument is crucial; however
the estimates for the error terms are more involved.
In \S5 we also bound a family of less singular positive operators which dominate operators 
in $\cI^{\r,-\s}$ when $\rho< 0$; thus we can then give a proof of  1.2.3.
Finally, in \S5, we discuss standard examples which show th sharpness of the results.
In \S 6 we establish  $L^p\to L^p$ bounds by suitable 
interpolation between  $L^2\to L^2$ and  Hardy-space estimates.
\S7 contains estimates for general operators in $\cI^{0,-\s}$ and 
additional interpolation arguments to finish the proof of Theorem 1.2.

\head{\bf 2. Preliminaries}
\endhead

\subheading{2.0. Notation}

2.0.1. $B_\eps$ will denote the open ball in $\Bbb R^d$ 
 of radius $\eps$ centered at the origin.

2.0.2. $m(D)$ denotes the convolution operator with Fourier multiplier $m(\xi)$. We split variables in
$\Bbb R^d=\bbR^{d-\ell}\times \bbR^\ell$ as $x=(x',x'')$ and denote by $h(D'')$ the 
convolution operator with Fourier multiplier $h(\xi'')$.

2.0.3. A function $F$ on $\{z:0\le \Re(z)\le 1\}$ is called 
 of admissible growth if $|F(z)|\le C e^{A|z|}$ for some $A>0$, $C\ge 0$. 

2.0.4.
The differentiability inequalities (1.11) are supposed to hold for all multiindices of length
$\le M_0$ where $M_0$ is large, say $M_0= 10^{100} d$, those multiindices are termed admissible.
Exponents
 $N, N_0,..., N_4$ in \S 4 and \S7 are assumed to be $\ge d+1$ and 
$\le  10^{10} d$.

2.0.5. 
We denote by $\zeta_0\equiv \om_0$ an even $C^\infty_0 (\Bbb R)$ function 
with  $\zeta_0(s)=1$ 
for $|s|\le 1/2$
and $\zeta_0(s)=0$ for $|s|\ge 1$.
Also let
$\zeta(s)=\zeta_0(s/2)-\zeta_0(s)$, $\om(s) 
=\om_0(s/4)-\om_0(s)$ 
so that $\zeta$ is supported in $[1/2,2]$ and   $\om $ is supported in  $[1/4,4]$; moreover
$$
\align &\zeta_0(s)+\sum_{j=1}^\infty\zeta(2^{-j}s)=1
\\
&\om_0(s)+\sum_{j=1}^\infty\om(4^{-j}s)=1
\endalign
$$ for all $s\in \Bbb R$.

2.0.6.
For two quantities $A$ and $B$ we write 
 $A\lc B$ or $B\gc A$ if there exists an absolute 
 positive constant $C$ so that $a\le C b$.
We write $A\approx B$ if both $A\lc B$  and $A\gc B$ hold.

\subheading{2.1 Standard assumptions}
For our Fourier integrals (1.10) and for the weakly singular Radon transforms
 any contribution away from the diagonal is handled by standard estimates for Fourier 
integral operators, see Lemma 3.1 below. Therefore,  in view of the compact support 
assumption on the kernel 
it is sufficient to prove
Theorems 1.1 and 1.2 under the assumption that the  kernels of our operators 
are supported in a  small neighborhood of a given point 
$(P,P)\in \Delta$.
We shall introduce coordinates that vanish at $P$, and assume that in these coordinates
the kernels are supported  where $|x|$, $|y|\le \eps^{10}$;
 $\eps$ is chosen in (2.16) below.

For further preparation  choose  $\eps_0>0$ so that in a neighborhood
of the closure of
$B_{\eps_0}\times B_{\eps_0}$  the  manifold $\cM$ is given as a graph
$$
y''=S(x,y');
\tag 2.1
$$
by performing a linear transformation we can also assume that
$$
S_{x'}(0,0)=O_{\ell,d-\ell}
\tag 2.2
$$
(the $\ell\times(d-\ell)$ zero-matrix).

Since $\Delta\subset \cM$ we have
$$
\align
x''=S(x,x') 
\tag 2.3
\endalign
$$
for all $x\in B_{\eps_0}$ and consequently
$$
\gather
S_{x'}(x,x')+S_{y'}(x,x')=0
\tag 2.4
\\
S_{x'x'}(x,x')+2S_{x'y'}(x,x')+S_{y'y'}(x,x')=0
\tag 2.5
\\
S_{x''}(x,x')= I_{\ell,\ell}
\tag 2.6
\endgather
$$
where
$ I_{\ell,\ell}$ denotes  the $\ell\times\ell$ identity matrix.

We shall also assume that for some constant $C_0\ge 1$
$$
\sum\Sb
|\alpha|\le 10^{100}d
\endSb
\sup\Sb|x|\le\eps_0 \\|y'|\le \eps_0\endSb 
|\partial_{x,y'}^\alpha S(x,y')|\le  
C_0.
\tag 2.7
$$
Moreover, by the assumption (1.9) and by (2.2) we have  for some positive 
 $c_0<1$
$$\|(\theta\!\cdot\! S_{x'y'}(0,0))^{-1}\|\le c_0^{-1},
\tag 2.8$$
for all unit vectors $\theta\in \bbS^{\ell-1}$; 
here $\|\cdot\|$ denotes the Hilbert-Schmidt norm.

\subheading{2.2 Straightening near the diagonal}

We now introduce a family of changes of variables,
 depending on unit vectors $u$ in $\Bbb R^{d-\ell}$
$$w\mapsto 
 \cQ(w;u):= (w',w''+F(w;u))$$
so that
$$\align
F_w'(0;u)&=0
\tag 2.9.1
\\
F(0;u)&= 0
\tag 2.9.2
\endalign
$$
and 
so that 
$$
\gathered
y''=S(x,y')\iff z''=
 \widetilde S(w,z'; u)
\\
\quad\text{ if  }\quad y=\cQ(z;u),
\quad x=\cQ(w;u)
\endgathered
\tag 2.10 $$
and
$$
\inn {u}{\nabla_{w'}} \widetilde S^{i}(w,w';u) =0, \quad i=d-\ell+1,\dots, d.
\tag 2.11
$$

To describe this change of variables let $B=B(u)$  be a rotation on 
$\Bbb R^{d-\ell}$ depending smoothly on $u$ such that
$Be_1=u$ (with  $e_1=(1,0,\dots,0)$).
We define an
 $\Bbb R^\ell$-valued function $G=G(\cdot\, ;u)$ by requiring that $G$ 
satisfies the following  
system of ordinary differential equations, with respect to the variable $w_1$
and initial  data depending on the parameters $w_2,\dots, w_d$:

$$
\gather
\frac{\partial G}{\partial w_1}(w)=
\inn{u}{S_{y'}}(Bw',w''+G(w),Bw') 
\\
G(0,w_2,\dots, w_d)= 0
\endgather
$$
Set 
$$F(w)\equiv F(w;u)=G(B^{-1}w',w'';u);
$$
then $F$ satisfies (2.9) and 
$$
 \inn{u}{\nabla_{w'}}F(w)=
\inn{u}{S_{y'}(w',w''+F(w),w')}. 
\tag 2.12 
$$


For the following discussion fix $u$.
Since the  functions $S$ and $\widetilde S$  are related by (2.10) we have
$$
\widetilde S(w,z')+F(z',\widetilde S(w,z))
= S(w',w''+F(w),z').
\tag 2.13
$$
Denote by $D_u=\inn{u}{\nabla_{z'}}$ the directional derivative with respect to $u$.
Differentiation of (2.13) yields
$$
D_u \widetilde S(w,z')+D_u F(z',\widetilde S(w,z))
+F_{z''}(z',\widetilde S(w,z)) D_u\widetilde S(w,z)=
 \inn{u}{\nabla_{y'}}S(w',w''+F(w),z')
$$
and by (2.12)  we obtain
$$\multline
D_u\widetilde S(w,z')+
\inn{u}{S_{y'}(z',\widetilde S(w,z')+F(z',\widetilde S(w,z')),z') }
+F_{z''}(z',\widetilde S(w,z')) D_u\widetilde S(w,z')
\\=
\inn{u}{S_{y'}(w',w''+F(w),z')}.
\endmultline
$$
Now we evaluate  for $w=z$ 
and take into account that $\widetilde S(z,z')=z''$  . This yields 
$$
(I+F_{z''}(z)) D_u\widetilde S(z,z')=0
$$
Since $F_{z''}(0)=0$ by (2.9), we obtain 
$\inn{u}{\nabla_{z'}}\widetilde S(z,z')=0$ in a neighborhood 
of $(0,0)$, and since also 
$\widetilde S_{w'}(w,w')+ \widetilde S_{z'}(w,w')=0$ this yields (2.11).

In view of (2.9) we may fix a number  $\delta_1\ll \eps_0$ so that
$$B_{\delta/2}\subset  \cQ(w,u) B_\delta \subset B_{2\delta} \quad \text{ for }
\quad \delta\le \delta_1, w\in B_\delta.
\tag 2.14$$

Let $$C_1=\sup_{|\alpha|\le 10^{100}d}\sup_{|w|\le \delta_1}|F^{(\alpha)}(w)|+C_0\tag 2.15$$
where $C_0$ is as in (2.7).
We may assume throughout this paper that the cutoff function $\chi$ in
(1.12) satisfies
$$\supp \chi\subset\{(x,y'): |x|+|y'|\le \eps\} \text{ where } 
0<\eps< (100d C_1/ c_0)^{-1}  \delta_1.
\tag 2.16$$
Moreover  the distribution kernels of the 
 the Fourier integrals 
defined by (1.10) are assumed to be supported  in $B_{\eps^{10}}\times B_{\eps^{10}}$.

Note also that for $|x|,|y|\le\eps$

$$
\gather
\|S_{x'}\|+\|S_{x''}-I_{\ell,\ell}\|
\ll \eps^{9} \ll \delta_1
\tag 2.17
\\
\big\|F_w\big\| \ll \eps^{9} \ll \delta_1.
\tag 2.18
\endgather
$$

\subheading{2.3. Adjoint operators} 
Suppose that  $\cM$ is given as a graph (1.7)with (1.8)  and the symbol has small
 $(x,y)$ support
then we may solve the equation $y''=S(x,y')$ in $x''$ so that
$y''=S(x', \fS(y,x'),y')$ and
$$y''-S(x',x'',y')= \cC(x,y)(x''-\fS(y',y'',x'))
\tag 2.19 $$ in a neighborhood of $\cM$,
with $\cC(x,y)$ is an invertible $\ell\times\ell$ matrix depending smoothly on $(x,y)$.
If  in the oscillatory integral (1.10) we make a linear change in the 
$\tau$-variables,
$\widetilde \tau= \cC(x,y)^T\tau$, then we see that (1.10) can be rewritten as a 
linear combination of integrals  with phase function $\inn{\tau}{x''-\fS(y,x')}$. 
This shows that 
for an operator in $\cI^{\rho,-\s}(\Om\times\Om,\cM,\Delta)$ the adjoint operator 
belongs to 
$\cI^{\rho,-\s}(\Om\times\Om,\cM^*,\Delta)$  where $\cM^*=\{(x,y):(y,x)\in\cM\}$
(and $\cM^*$ satisfies (1.3), (1.4)).

\head{\bf 3. Nonsingular Radon transforms and scaling}\endhead 

We first recall a well-known result on $L^p\to L^q$ estimates for 
Fourier integral operators
associated to a canonical graph.
These estimates take care 
of  contributions of the kernels away from the diagonal.
In the formulation of this Lemma the order of a Fourier integral opertator 
is as in the 
standard theory of Fourier integral operators; 
thus the standard Radon-type operators is of order 
$-(d-\ell)/2$.

\proclaim{Lemma 3.1 }
 Suppose $-\ell< \r <\frac{d-\ell}2$.

 Let $T$ be a Fourier integral operator of order
$\r-\frac{d-\ell}2$ associated 
to a local canonical graph $\Cal C\subset
T^*\Omega\setminus \{0\}\times 
T^*\Omega\setminus \{0\}$. 
Suppose that the restrictions  $\Cal C$ of  the
 projections $(x,y)\to x$ and $(x,y)\to y$
have differentials with maximal rank $d$ and that the projection
$\Cal C\to \Omega\times\Omega$ has a differential with constant 
 rank $\le 2d-\ell$.
Suppose that the distribution kernel of $T$ has compact support.

(i) If  $\rho>0$ then
$T$ maps  $L^p$ to $L^q$ 
if
 $(1/p, 1/q)$ belongs to the closed triangle with corners
$(\frac{\r}{d-\ell}, \frac{\r}{d-\ell})$,
$(\frac{d-\ell-\r}{d-\ell},\frac{d-\ell-\r}{d-\ell})$ and 
$(\frac{d-\r}{d+\ell}, \frac{\r+\ell}{d+\ell})$.

(ii) If $\rho=0$ then 
$T$ maps  $L^p$ to $L^q$ 
if
 $(1/p, 1/q)$ belongs to the closed triangle with corners
$(0,0)$, $(1,1)$ and 
$(\frac{d}{d+\ell}, \frac{\ell}{d+\ell})$, with the possible exception of
the corners $(0,0)$ and $(1,1)$; 
then an $H^1\to H^1$ or $L^\infty\to BMO$ bound holds.

(iii) 
If $-\ell<\rho<0$ then  $T$ maps $L^p$ to $L^q$ if 
 $(1/p, 1/q)$ belongs to the pentagon with  corners
$(1,1)$, $(0,0)$, $(1,\frac{\ell+\rho}{\ell})$, $(\frac{-\rho}{\ell},0)$ and
$(\frac{d-\r}{d+\ell}, \frac{\r+\ell}{d+\ell})$, with the possible exceptions of the points
$(1,\frac{\ell+\rho}{\rho})$, $(\frac{-\rho}{\ell},0)$.
\endproclaim

\demo{Sketch of the argument}
The main $L^p\to L^{p'}$ estimates are essentially proved in 
\cite{1}.
We sketch the argument. Consider first
the main  endpoint $L^{\frac{d+\ell}{d-\rho}}\to 
L^{\frac{d+\ell}{\r+\ell}}$ estimate. In view of the constant rank assumptions on the 
projection of $\cC$  to the base space we may after appropriate localization  and choice of coordinates 
write the 
kernel  as the sum  $\sum_{k\ge 1} K_k(x,y)$ and a $C^\infty_0$ function; here
$$
K_k(x,y)=2^{k\rho} \int e^{i\inn{\tau}{y''-S(x,y')}} a_k(x,y,\tau) d\tau
$$
where the integral is extended over a conic open set of $\Bbb R^\ell$, $S$ is as in the
 introduction,  the symbols $a_k$ are of  order $0$ with uniform  bounds in 
$k\ge 1$, and $a_k(x,y,\theta)=0$ if $|\theta|\notin (2^{k-1}, 2^{k+1})$.

Let $T_k$ be the operator with kernel $K_k$. 
Standard $L^2$ theory (see \cite{12}, \cite{22}) shows that
$T_k$ is bounded on  $L^2$, with norm $O(2^{k(\rho-\frac{d-\ell}2)})$.
Clearly $|K_k(x,y)|\lc 2^{k(\rho+\ell)}$.
Thus $T_k$ maps $L^1$ to $L^\infty$ with norm
$\lc  2^{k(\rho+\ell)}$.
Interpolation yields that
$T_k$
maps
 $L^{\frac{d+\ell}{d-\rho}}$ to 
$L^{\frac{d+\ell}{\r+\ell}}$ with bounds uniform in $k$. 
Since we assume that the canonical relation $\cC$ does not meet $\{0\}\times T^*\Om$ and
$T^*\Om\times\{0\}$ one can use standard integration by parts arguments 
(\cite{12}) and 
Littlewood-Paley theory 
to put the pieces together
and one obtains the desired 
 $L^{\frac{d+\ell}{d-\rho}}\to 
L^{\frac{d+\ell}{\r+\ell}}$ estimate, {\it cf.}  also \cite{1}. 
For the endpoint $L^p\to L^p$ (or $H^1\to L^1$ estimate)
and more references see \cite{22, ch. IX}.

Finally assume $-\ell<\rho<0$. Then an integration by parts argument 
shows that
$$|K_k(x,y)|\lc 2^{k(\ell+\rho)}(1+2^k|y''-S(x,y')|)^{-N}$$
and therefore the sum in $k$ is bounded by $|y''-S(x,y')|^{-(\rho+\ell)}$.
In view of the compact support of the kernel we see that
$K(x,\cdot)$ and $K(\cdot, y)$ are uniformly in Weak-$L^{\frac{\ell}{\ell+\rho}}$.
Thus the operator maps $L^1$ to Weak-$L^{\frac{\ell}{\ell+\rho}}$. A similar argument applies to the adjoint operator.
Now one uses the Marcinkiewicz interpolation to interpolate with the endpoint
 $L^{\frac{d+\ell}{d-\rho}}\to 
L^{\frac{\r+\ell}{d+\ell}}$ estimate and further interpolation with the trivial $L^1$ and $L^\infty$ estimates to conclude.\qed
\enddemo

Let $\chi\in C^\infty_0(\bbR^d\times\bbR^d)$ be a nonnegative function.
Now let $-\ell<\rho\le 0$, $0<\s<d-\ell$. If also $\rho<0$ we  define 
the distribution kernel  $G^{\rho,\sigma}$ by
$$G^{\rho,\s}(x,y) =  \chi(x,y)
|x'-y'|^{-(d-\ell-\s)}
 c_{\ell,\rho}|y''-S(x,y')|^{-(\rho+\ell)}
\qquad \text{ if } -\ell<\rho<0
\tag 3.1
$$
where
$$c_{\ell,\rho}=2^\rho\pi^{-\ell/2} \frac{\Gamma(\tfrac{\ell+\r}2)}
{\Gamma(\tfrac{-\r}2)}
$$
so that the Fourier transform on $\bbR^\ell$ of 
$c_{\ell,\rho}|\cdot|^{-(\ell+\r)}$ is 
$|\xi|^{\r}$, see \cite{8}. Define 
$G^{0,\s}=\lim_{\rho\to 0-} G^{\rho,\s}$
where the limit is taken in the sense of distributions; clearly
$$G^{0,\s}(x,y) =
\delta(y''-S(x,y'))
|x'-y'|^{-(d-\ell-\s)}  \chi(x,y) .\tag 3.2
$$
Define  the operator $\cR^{\r,\s}$ by
$$\cR^{\r,\s} f(x)= \inn {G^{\rho,\s}(x,\cdot)}{f}
\tag 3.3
$$ so that for $\rho=0$ we recover the weakly singular Radon transform.
We wish  to  apply 
Lemma 3.1 to dyadic pieces localized in $x'-y'$,
after a suitable rescaling.
Therefore we decompose dyadically 
$$\cR^{\r,\s}=\sum_j (\cR_j^{\r,\s}+\cE_j^{\rho,\s})
\tag 3.4$$
with
$$\align
\cR_j^{\rho,\s} f(x)&=
2^{j(d-\ell-\sigma)}
\inn {b_j^{\rho,\s}(x,\cdot)}{f}
\tag 3.4
\\
\cE_j^{\rho,\s} f(x)&=
2^{j(d-\ell-\sigma)}
\inn {h_j^{\rho,\s}(x,\cdot)}{f}
\tag 3.5
\endalign
$$
where
$$
b_j^{\r,\s}(x,y)=
2^{-j(d-\ell-\s)} \zeta(2^j|x'-y'|)
\zeta_0(\tfrac{|y''-S(x,y')|}{|x'-y'|^2}) G^{\rho,\s}(x,y)
\tag 3.6$$
and
$$
h_j^{\r,\s}(x,y)=
2^{-j(d-\ell-\s)} \zeta(2^j|x'-y'|)
(1-\zeta_0(\tfrac{|y''-S(x,y')|}{|x'-y'|^2})) G^{\rho,\s}(x,y).
\tag 3.7
$$
Note that this implies $h^{0,\s}\equiv 0$.

\proclaim{Proposition 3.2} 
Let $0<\s<d-\ell$, $-\ell<\r\le 0$  and let $\cR_j^{\r,\s}$ be as in (3.1).

(i) Suppose that
 $(1/p, 1/q)$ belongs to the triangle with corners
$(0,0)$,
$(1,1)$ and 
$(\frac{d}{d+\ell}, \frac{\ell}{d+\ell})$.
Then 
$$\|\cR_j^{0,\s} f\|_q\lc 2^{j[(d+\ell)(\frac 1p-\frac 1q)-\sigma]}\|f\|_p. 
$$

(ii) 
Suppose that $-\ell <\r<0$. Then the inequality
$$\|\cR_j^{\rho,\s} f\|_q\lc 2^{j[(d+\ell)(\frac 1p-\frac 1q)+2\rho-\sigma]}\|f\|_p 
$$
holds if  
 $(1/p, 1/q)$ belongs to the pentagon with  corners
$(1,1)$, $(0,0)$, $(1,\frac{\ell+\rho}{\ell})$, $(\frac{-\rho}{\ell},0)$ and
$(\frac{d-\r}{d+\ell}, \frac{\r+\ell}{d+\ell})$, with the possible exception of the points
$(1,\frac{\ell+\rho}{\ell})$, $(\frac{-\rho}{\ell},0)$.
Then

\endproclaim

\demo{\bf Proof}
Let $\delta>0$ and $$B(a,\delta)=\{y: |y'-a'|\le \delta, 
|y''-a''-\inn{S_{y'}(a,a')}{y'-a'}|\le \delta^2\}.
\tag 3.8$$
A sufficiently small
 neighborhood $U$ of the origin is then made into a space of homogeneous
space with the balls $B(x,\delta)$ (see \cite{16}, \cite{22} at least for the case 
$\ell=1$), and
for sufficiently large $j$ we can cover $U$  with a family of balls
 $B(x_\nu, 2^{-j})$ which have bounded overlap.

Fix $j$ and observe that if $f$ is supported in $B(x_\nu, 2^{-j})$ then
$\cR_j f$ is supported in $B(x_\nu, C 2^{-j})$ for a fixed $C$.
Therefore in order to prove the asserted inequality it suffices to verify 
it under the assumption that $f$ is supported in a ball
 $B(a,\delta)$ where $a\in \Omega$ is near the origin.

Fix $a$. Then we perform  an affine  change of variables, 
so that
in the new coordinates  we can write $R_j^{\r,\s}$ as
in (3.4), (3.6) with $S(x,y_1) $ replaced by $s(x,y_1)$ satisfying
$$s_{x'}(a,a')=0, \qquad 
s_{y'}(a,a')=0.
\tag 3.9
$$
(3.9) implies that  the ball $B(a,2^{-j})$ is contained in
$$\{y:|y'-a'|\le A 2^{-j} ,|y''-a''|\le A 2^{-2j}\}$$ for suitable $A$.
Moreover we  also  see
 the rotational curvature in (1.9)  at $(a,a')$ 
is given by  $\det \theta\cdot s_{x'y'}(a,a')$ 
since we still have $s_{x''}(a,a')=I_{\ell,\ell}$, {\it cf.} (2.6).

We now perform a scaling argument and
write
$$
\cR_j^{\r,\s} f(a'+2^{-j}v',a''+ 2^{-2j}v'')= 2^{j(2\r- \s)} 
\widetilde \cR_j^{\r,\s} f_j(v)
$$
where
$$
\align
\widetilde\cR_j^{\r,\s} g(v)&= \inn{\widetilde b_j^{\r,\s}(v,\cdot)}{g},
\\
f_j(w', w'')&= f(a'+2^{-j}w', a''+2^{-2j} w''),
\\
S_{j,a}
(v,w')&= 2^{2j}(-a''+ s(a'+2^{-j}v', a''+2^{-2j}v'', a'+2^{-j} w')),
\\
\widetilde b_j^{\r,\s}(v,w) &= 
b_j^{\r,\s}(a'+2^{-j}v', a''+2^{-2j} v'', a'+2^{-j}w', a''+2^{-2j} w'').
\endalign
$$

In view of $s(a,a')=a''$ and 
$s_{x'}(a,a')=0$ we check that the derivatives of $S_{j,a}$ are uniformly bounded 
(in a fixed neighborhood of $(0,0)$, which can be chosen independently of $j$ and $a$)
 and also that the 
rotational curvature is bounded
below.

The rescaled operators $\widetilde\cR_j^{\r,\s}$ are standard Fourier integral operators, to which 
Lemma 3.1 (ii), (iii) can be applied, the resulting $L^p\to L^q$ bounds are uniform in $j$,
and in $a$.
We apply  
Lemma 3.1 with the relevant choice of $p$ and $q$ and it follows that
$$2^{j\frac{d+\ell}q}\|\cR_j^{\r,\s} f\|_q
= 2^{j(2\r-\s)}\|\widetilde\cR_j f_j\|_q\lc
 2^{j(2\r-\s)}\|f_j\|_p\lc
2^{j(2\r-\s)}2^{j\frac{d+\ell}p}\|f\|_p
$$
which proves the Proposition. \qed
\enddemo

For the estimation of the error term involving the terms  $\cE_j^{\rho,\s}$ see Proposition 4.2 
below.

\head{\bf 4. 
Regular and product type fractional integrals}\endhead
 
In this section we study nonisotropic and  product type pseudodifferential operators, which
come up as low frequency contributions to operators in $\cI^{\r,-\s}$; in particular we prove 
$L^p\to L^q$ estimates for the error term in (3.4).
We recall a sharp version  of Young's inequality  (see Theorem (6.35) in \cite{6}) which states that
the conditions $1<p<q<\infty$ and 
$$
\sup_x\|K(x,\cdot)\|_{L^{r,\infty}}
+
\sup_y\|K(\cdot,y)\|_{L^{r,\infty}}<\infty, \qquad 
\frac 1r= 1-\frac 1p+\frac 1q,
\tag 4.1
$$
imply that the integral operator with kernel $K(x,y)$ is bounded from $L^p\to L^q$.

\proclaim{Lemma 4.1} Suppose $1<p<q<\infty.$ Define
$$
\align
K^{\r,\s}_1(x,y)&= \chi(x,y)|x'-y'|^{\sigma-d+\ell}|y''-S(x,y')|^{-\r-\ell}
\tag 4.2
\\
K^{\r,\s}_2(x,y)&= \chi(x,y)(|x'-y'|+|y''-S(x,y')|^{1/2})^{\s-2\r-d-\ell}
\tag 4.3
\endalign
$$
and
$$
K^{\r,\s}_3(x,y)=\cases
 \chi(x,y)|x'-y'|^{\sigma-d+\ell}|y''-S(x,y')|^{-\rho-\ell}
&\text{ if } |x'-y'|^2\le 10|y''-S(x,y')|
\\
0&\text{ if } |x'-y'|^2\ge 10|y''-S(x,y')|
\endcases .
\tag 4.4
$$

(i) Assume  $0<\sigma<d-\ell$, $-\ell<\rho<0$, $(d-\ell)(1/p-1/q)\le\sigma$, 
$\ell(1/p-1/q)\le -\r$.
Then the integral operator with kernel $K_{1}^{\r,\s}$ maps $L^p$ to $L^q$.

(ii) 
 Assume   $-\ell<\rho\le 0$, $0<\sigma<d-\ell$ 
 and  $(d+\ell)(1/p-1/q)\le\sigma-2\r$.
Then
 the integral operator with kernel $K_{2}^{\r,\s}$  
maps $L^p$ to $L^q$.

(iii)  Assume 
$-\ell<\rho\le 0$,
  $-\r(d-\ell)/\ell<\sigma<d-\ell$ and   $(d+\ell)(1/p-1/q)\le\sigma-2\r$.
Then
 the integral operator with kernel $K_{3}^{\r,\s}$  
maps $L^p$ to $L^q$.

\endproclaim

\demo{\bf Proof}
We first consider  (i). Let $\cJ_{x',y'} $ 
denote the integral operator acting on functions in $\bbR^\ell$, 
with kernel
$$J_{x',y'}(x'',y'')= \chi(x',x'',y',y'') |y''-S(x,y')|^{-\r-\ell}.$$
If $\ell(1/p-1/q)\le -\rho$ then $\sup_{x''}\|J_{x',y'}(x'',\cdot)\|_{L^{r,\infty}}\le C$ 
for $1/r=1-1/p+1/q$, uniformly in $x',y'$. Since 
the quantities
 $|y''-S(x,y')|$ and $|x''-\fS(y,x')|$ are comparable ({\it cf.} \S2.3) 
we also have 
$\sup_{y''}\|J_{x',y'}(\cdot,y'')\|_{L^{r,\infty}}\le C$. Thus by the  
sharp form of Young's inequality stated above   the condition
$\ell(1/p-1/q)\le -\r$ implies that 
$J_{x',y'}$ maps $L^p(\bbR^\ell)$ to $L^q(\bbR^\ell)$, with
bounds independent of $x',y'$. 
Likewise, since $(d-\ell)(1/p-1/q)\le \s$ the integral operator 
with kernel
$\widetilde \chi(x',y') |x'-y'|^{\s-d+\ell}$ maps $L^p(\Bbb R^{d-\ell})$ to 
$L^q(\bbR^{d-\ell})$ 
 if $\widetilde \chi$ is compactly supported. Thus
by Minkowski's inequality (if $T_1^{\r,\s}$ is the integral operator with kernel $K_1^{\r,\s}$)
$$
\align
\|T_1^{\r,\s}f\|_q&\le\Big(\int\Big [\int 
\widetilde \chi(x',y') |x'-y'|^{\s-d+\ell} \big\|J_{x',y'}[f(y',\cdot)]
\big\|_{L^q(\Bbb R^\ell)} dy' 
\Big]^q dx'\Big)^{1/q}
\\
&\lc\Big(\int\Big [\int 
\widetilde \chi(x',y') |x'-y'|^{\s-d+\ell} \|f(y',\cdot)\|_{L^p(\Bbb R^\ell)}
 dy' 
\Big]^q dx'\Big)^{1/q}
\\
&\lc\Big(\int \|f(x',\cdot)\|_{L^p(\Bbb R^\ell)}^p dx' 
\Big)^{1/p}
\endalign$$
and hence $T_1^{\r,\s}$ is bounded from $L^p(\bbR^d)$ to $L^q(\bbR^d)$. This proves (i).

(ii) is proved by checking directly 
the condition (4.1) for $r\le \tfrac{d+\ell}{d+\ell+2\r-\s}$;
the calculation is standard and therefore omitted.

It remains to consider the operator with kernel $K^{\rho,\sigma}_3$.
We now  fix $x$ and prove
$\|K^{\rho,\sigma}_3(x,\cdot)\|_{L^{r,\infty}}\le C$ with $C$ independent of $x$; here again
$r= \tfrac{d+\ell}{d+\ell+2\r-\s}$.
 Let $v'=y'-x'$ and $v''=y''-S(x,y')$.

For $\alpha>0$ let 
$$\Om(\alpha)=\{(v',v''):
|v'|^{\sigma-d+\ell}|v''|^{-\r-\ell}>\alpha,   \ |v'|^2\le 10|v''|^{1/2}, \ |v''|\le C_2\}.
$$ 
We have to show that the  set $\Om(\alpha)$ 
has measure $O(\alpha^{-r})$. 
If $v\in \Omega(\alpha)$ then $|v'|^2\le 10|v''|\le 10 
\alpha^{-\frac{1}{\r+\ell}}
|v'|^{-\frac{d-\ell-\sigma}{\rho+\ell}}
$ and this  implies 
$|v'|^{2+\frac{d-\ell-\s}{\ell+\r}}\lc \alpha^{-\frac{1}{\rho+\ell}}$ or
$|v'|\lc \alpha^{-\frac{1}{d+\ell+2\r-\s}}$.
Thus
$$
|\Om(\alpha)|\lc 
\int\limits_{|v'|\lc \alpha^{-\frac{1}{d+\ell+2\r-\s}}}
\alpha^{-\frac{\ell}{\r+\ell}}|v'|^{-\frac{d-\ell-\sigma}{\rho+\ell}\ell} dv'
$$
Now the condition
$-\r\frac{d-\ell}{\ell}<\sigma$ is equivalent with $-\frac{d-\ell-\s}{\rho+\ell}\ell> -(d-\ell)$
and therefore one can verify
$$
|\Om(\alpha)|\le C \alpha^{-\frac\ell{\ell+\r}} 
\alpha^{-\frac{1}{d+\ell+2\r-\s}(d-\ell-\frac{(d-\ell-\s)\ell}{\ell+\r})}=C
\alpha^{-\frac{d+\ell}{d+\ell+2\r-\s}}
$$
and thus
$\sup_x\|K^{\rho,\sigma}_3(x,\cdot)\|_{L^{r,\infty}}<\infty$.
The verification of the condition 
$\sup_y\|K^{\rho,\sigma}_3(\cdot,y)\|_{L^{r,\infty}}<\infty$ is similar.
\qed
\enddemo

\proclaim{Proposition 4.2}
Suppose that $1< p\le q< \infty$,  $-\ell <\r<0$ and  $0<\s<d-\ell$.
Let $\cE^{\r,\s}=\sum_j \cE^{\r,\s}_j$ (as defined in (3.5))
Then $\cE^{\r,\s}$ is bounded from $L^p$ to $L^q$ if either one of the following two 
conditions  is satisfied.

(i) $-\rho\frac{d-\ell}\ell< \s<d-\ell$ and $(d+\ell)(1/p-1/q)\le \s-2\r$.

(ii) $0<\s\le -\rho\frac{d-\ell}\ell$ and $(d-\ell)(1/p-1/q)\le \s$.
\endproclaim

\demo{\bf Proof}
The kernel  of $\cE^{\r,\s}$ can be estimated by both $K_1^{\r,\s}$ and
$K_3^{\r,\s}$ in Lemma 4.1. For (i) apply the estimate for the integral operator with kernel
$K_3^{\r,\s}$. To prove (ii) from Lemma 4.1  observe that 
inequality  $\ell(1/p-1/q)\le -\r$  is implied 
by $0<\s\le -\rho\frac{d-\ell}\ell$ and $(d-\ell)(1/p-1/q)\le \s$.\qed
\enddemo

We shall now look at the basic dyadic pieces in decompositions of operators in $\cI^{\r,-\s}$.
Let
$$
\beta_{k,m}(x,y,\tau,\xi) =
\cases
 \om(2^{-2k}|\tau|)
 \zeta(2^{-m}|\xi|) 
&\quad\text{ if } k>0, m>0,
\\ \om_0(|\tau|)
 \zeta(2^{-m}|\xi|) 
&\quad\text{ if } m>0, k=0,
\\
 \om(2^{-2k}|\tau|)
 \zeta_0(|\xi|) 
&\quad\text{ if } k>0, m=0,
\\
 \om_0(|\tau|)
 \zeta_0(|\xi|) 
&\quad\text{ if } k=m=0.
\endcases
\tag 4.5$$
Suppose $a\in S^{\rho,-\sigma}$. Let 
$$
K_{k,m}(x,y)=\iint_{\Bbb R^{d-\ell}\times \Bbb R^\ell}
 e^{\imath[\inn{\tau}{y''-S(x,y')}+\inn{\xi}{x'-y'}]}
(a \beta_{k,m})(x,y,\tau,\xi) d\tau d\xi.
\tag  4.6
$$
Let $T_{k,m}$ be the integral operator with kernel $K_{k,m}(x,y)$.

\proclaim{Lemma 4.3} If $a\in S^{\rho,-\sigma}$ then

(i) $$
|K_{k,m}(x,y)|
\lc 2^{2k\rho-m\sigma} \frac {2^{2k\ell}}{(1+ 2^{2k}|y''-S(x,y')|)^{N}}
\frac {2^{m(d-\ell)}}{(1+2^m|x'-y'|)^N};
\tag 4.7
$$

moreover
$$|\nabla K_{k,m}(x,y)|\lc \max\{ 4^{k}, 2^{m}\}
 2^{2k\rho-m\sigma} \frac {2^{2k\ell}}{(1+ 2^{2k}|y''-S(x,y')|)^{N}}
\frac {2^{m(d-\ell)}}{(1+2^m|x'-y'|)^N}.
\tag 4.8
$$

(ii) Let $K$ be the Schwartz kernel of an operator in $\cI^{\rho,-\sigma}$ given by (1.10),
 and assume that $-\ell<\rho<0$, $0<\sigma<d-\ell$. 
Then $K$ satisfies 
$$
|K(x,y)|\lc  |y''-S(x,y')|^{-\r-\ell} 
|x'-y'|^{\s-d+\ell}.
$$
\endproclaim
\demo{\bf Proof} (i) follows by integration by parts. (ii) is deduced from (i)
by summing geometric series.\qed
\enddemo

\subheading{Proof of Theorem 1.2.4} Immediate from Lemma 4.3 (ii) and Lemma 4.1.\qed

We shall now look at a general operator in $\cI^{\r,-\s}$ and consider the 
contribution  which gives rise to a nonisotropic pseudo-differential operator .


\proclaim{Proposition 4.4} 
Let $a\in S^{\rho,-\sigma}$ and suppose that
 $1<p\le q<\infty$.
Suppose that $-\ell<\rho\le 0$
and that $-\rho\frac{d-\ell}\ell< \s<d-\ell$ and $(d+\ell)(1/p-1/q)\le \s-2\r$.
%
%
Then the operator
$\sum_{k\ge 0}\sum_{m\ge k} T_{k,m}$ is bounded from $L^p$ to $L^q$.
\endproclaim

\demo{Proof} We use the kernel estimates (4.7)  and sum. We find that the kernel $P(x,y)$ of 
$\sum_{k\ge 0}\sum_{m\ge k} T_{k,m}$ satisfies the estimate
$$
|P(x,y)|\lc\cases
|x'-y'|^{\s-2\rho-d-\ell}\quad&\text{ if } |y''-S(x,y)|^{1/2}\lc |x'-y'|
\\
|x'-y'|^{\ell-d+\s}|y''-S(x,y')|^{-\rho-\ell}\quad&\text{ if } |y''-S(x,y)|^{1/2}\gc |x'-y'|
\endcases .
$$
Thus
$$|P(x,y)|\lc K_2^{\rho,\s}(x,y)
+K_3^{\rho,\s}(x,y)$$ 
and the assertion   follows from Lemma 4.1.  \qed
\enddemo

For later use we also write down a similar estimate for an operator with localization 
in $|x'-y'|$.

\proclaim{Lemma 4.5} 
Let $a\in S^{0,-\sigma}$ and $K_{k,m}$ as in (4.6), with $\rho=0$. 
Denote by 
$W_{k,m}$  the operator with kernel $K_{k,m}(x,y)\zeta_0(2^{k}(|x'-y'|))$.
Suppose  $1<p\le q<\infty$ and
$(d+\ell)(1/p-1/q)\le \sigma$, $0<\sigma<d-\ell$. 
Then for $s>0$  the operator
$\sum_{k>s} W_{k,k-s}$ is bounded from
from $L^p$ to $L^q$, with operator norm $O(2^{-s (d-\ell- \s)})$.
\endproclaim

\demo{\bf Proof} This follows in a  straightforward manner  from (4.7) and
 Lemma 4.1. We have the estimate
$$
|K_{k,k-s}(x,y)|\lc\frac{2^{k(d+\ell-\s)} 2^{-s(d-\ell-\s)}}
{(1+2^{2k}|y''-S(x,y')|+2^{k-s}|x'-y'|)^N}|\zeta_0(2^k(x'-y'))|;
$$
here we choose $N>(d+\ell-\s)$.
If  $|y''-S(x,y')|\le |x'-y'|\approx 2^{-k}$ we simply dominate by
$2^{k(d+\ell-\s)} 2^{-s(d-\ell-\s)}$ which is in the present case controlled  by 
$2^{-s(d-\ell-\s)} K_2^{\rho,\sigma}(x,y)$ ({\it cf.} (4.3)).

If $|y''-S(x,y')|\ge |x'-y'|\approx 2^{-k}$ then
$|K_{k,k-s}(x,y)|\lc 
2^{-s(d-\ell-\s)}|y''-S(x,y')|^{-(d+\ell-\s)/2}$ and 
in the case under consideration this is also 
controlled  by
$2^{-s(d-\ell-\s)} K_2^{\rho,\sigma}(x,y)$. Since for fixed $(x,y)$ the sum
$\sum_{k>s} K_{k,k-s}(x,y)$ contains at most three terms, we see that the assertion follows from Lemma 4.1.
\qed
\enddemo

\head{\bf 5. Weakly singular Radon transforms and some variants}\endhead

In this section we give a proof of Theorem 1 and part 1.2.3 of Theorem 1.2.
We first introduce an additional angular localization in the angular variable.

Let  $v\in \Bbb R^{d-\ell}$ be a unit vector.
 Let
$$\align
 \kappa(x,y)&=
\chi(x,y)
\zeta_0\big(\eps^{-10}\big|\frac{x'-y'}{|x'-y'|}-v\big|\big)
\zeta_0\big(\frac{|y''-S(x,y')|}{|x'-y'|^2}\big)
\\
\kappa_j(x,y)&= \kappa(x,y)\zeta(2^j(|x'-y'|)
\endalign
$$
here $\chi$ is  a 
nonnegative smooth function supported where $|x|+|y'|\le \eps^{10}$ (see (2.16)).
Thus 
$$\supp    \kappa
\subset\big\{(x,y): \big|\frac{x'-y'}{|x'-y'|}-v\big| 
\ll \eps^{10}, |x|\le \eps^{10}, |y|\le\eps^{10}, |y''-S(x,y')|\le |x'-y'|^2\big \}.
\tag 5.1
$$

Let $G^{\r,\s}$ be as in (3.1) and define
$$
R^{\r,\s} f(x)= \inn{G^{\r,\s}(x,\cdot) \ka(x,\cdot)}{f}.
\tag 5.2
$$

The operator $\cR^{0,\sigma}$ introduced in \S3
is a finite sum of operators of type $R^{0,\sigma}$ (with suitable choices of $\chi$ and $v$). Moreover, for $\rho<0$ we recover the operators 
$\cR^{\r,\s}$ modulo error terms which are  already estimated by 
Proposition 4.2. The case $\rho=0$ of the following  result implies the
assertion of  Theorem 1.1.

\proclaim{Theorem 5.1} Let $1\le p\le q\le \infty$.

(i)  Suppose that  
 $(1/p, 1/q)$ belongs to the intersection of the  halfspace defined by
 $(d+\ell)(\frac 1p-\frac 1q)\le 
\s$ with the triangle with
corners $(0,0)$, $(1,1)$ and 
$(\frac{d}{d+\ell}, \frac{\ell}{d+\ell})$. Then 
$R^{0,\s}$ maps $L^p$ to  $L^q$.

(ii)  Suppose $-\ell< \r< 0$ and  $-\r\frac{d-\ell}{\ell}< \s<d-\ell $.
 Suppose that  
 $(1/p, 1/q)$ belongs to the intersection of the  halfspace defined by
 $(d+\ell)(\frac 1p-\frac 1q)\le 
\s-2\r$ with the  pentagon with  corners
$(1,1)$, $(0,0)$, $(1,\frac{\ell+\rho}{\ell})$, $(\frac{-\rho}{\ell},0)$ and
$(\frac{d-\r}{d+\ell}, \frac{\r+\ell}{d+\ell})$, with the  exception of the points
$(1,\frac{\ell+\rho}{\ell})$, $(\frac{-\rho}{\ell},0)$.
Then $R^{\r,\s}$ maps $L^p$ to  $L^q$.
\endproclaim

For the rest of this section we fix $\rho, \sigma$ and will not explicitely indicate 
the dependence on these parameters. If  $p=q$ the assertion is easily verified 
 by Minkowski's inequality. This also applies to the cases
  $p=1$ and $q<\ell/(\ell+\rho)$, and
$q=\infty$ and $p<-\ell/\rho$  (when $-\ell<\r<0$). Thus we may assume 
 $1<p<q<\infty$,
and
that
$(1/p,1/q)$ satisfies the restrictions in Theorem 5.1;
 moreover  we may assume 
 $p\le 2$ since the case  $p>2$ follows by considering 
the adjoint operator. It is always assumed that the function $f$ is supported where $|y|\le \eps^{10}$ and $\eps $ is as in (2.16). These assumptions are always assumed but not explicitly 
stated in various lemmas throughout this section.

Define
$$\cR_j f(x)=\inn {G^{\r,\s}(x,\cdot)\kappa_j(x,\cdot)}{f}.
\tag 5.3$$
Then   $\cR_j$ is bounded from 
$L^p$ to $L^q$ with a bound independent of $j$, by Proposition 3.2. Let $M$ be such that
$2^M\ge (\eps c_0)^{-10} $ (with $c_0$ as in (2.8)) and let $J$ be a finite set of  integers, all of them $\ge M$. Let 
$$
\cR f=\sum_{j\in J}\cR_j f.
\tag 5.4
$$
A priori we know that $\cR$ is bounded from $L^p\to L^q$ with norm
 $O(\card(J))$, and our task is to improve this to show   that the $L^p\to L^q$ bound
is independent of the cardinality of $J$. 
Once this is proved the 
$L^p\to L^q$ boundedness of $\cR^{\rho,\sigma}$ follows 
immediately from applications
of  the monotone convergence theorem.

We begin by cutting out the low frequencies 
(here we follow essentially \cite{2}, \cite{11}) and split $\cR=\cA+\cB$ with
$$\align
\cA&=
\sum_{j\in J} \omega_0(2^{-2j}|D''|)\cR_j,
\tag 5.5.1
\\
 \cB&=
\sum_{j\in J} (I- \omega_0(2^{-2j}|D''|))\cR_j.
\tag 5.5.2
\endalign
$$

We  first prove
\proclaim{Lemma 5.2} The operator 
$\cA$ 
is bounded from 
$L^p$ to $L^q$, with norm independent of the family $J$.
\endproclaim

\demo{\bf Proof}
Since the convolution kernel 
$\omega_0(2^{-2j}|D''|)$ is $O(2^{2j\ell}(1+2^{2j}|x''|)^{-N}$ we see that  for $\rho<0$
$$
\align
&\big|\omega_0(2^{-2j}|D''|)\cR_j f(x)\big|
\\&\lc 
\iiint\limits\Sb|y''-S(x',w'',y')|\lc 2^{-2j}\\|x'-y'|\approx 2^{-j}\endSb
 \frac{2^{2j\ell}}{(1+2^{2j}|x''-w''|)^{N}} 
G^{\r,\s}(x',w'',y',y'') |f(y',y'')| dy' dy'' dw''
\\
&\lc 
\iint\limits\Sb\{(y',y''):\\|x'-y'|\approx 2^{-j}\}\endSb
 |f(y',y'')| 2^{j(d-\ell-\s)} 
\int\limits\Sb|w''-\fS(y',y'',x')|\\ \lc 2^{-2j}\endSb
 \frac{2^{2j\ell}}{(1+2^{2j}|x''-w''|)^{N}} 
|w''-\fS(y',y'',x')|^{-\rho-\ell} dw'' \, dy' dy'' 
\\
&\lc 
\iint\limits\Sb|x'-y'|\approx 2^{-j}\endSb
 \frac{2^{j(d+\ell+2\r-\sigma)}}{(1+2^{2j}|x''-\fS(y',y'',x')|)^{N}}
|f(y',y'')| dy' dy'';
\endalign
$$
here $\fS$ is as in \S2.3.
The same estimate applies to the case $\rho=0$ (with only notational changes in the argument).

We see that
the  kernel of 
$\omega_0(2^{-2j}D'')\cR_j$
can be estimated by $K_2^{\r,\s}$ (as in (4.3)), uniformly in $j$. This bound also applies to the sum
$\sum_{j\in J}\omega_0(2^{-2j}D'')\cR_j$ since the kernel of 
$\omega_0(2^{-2j}D'')\cR_j$ is supported where $|x'-y'|\approx 2^{-j}$.
 Thus the assertion follows from Lemma 4.1.
\qed
\enddemo

We now turn to the operator $\cB$ and we  shall first prove 
estimates for a  frequency localized variant.

\proclaim{Proposition 5.3}
Let $\vth$ be a fixed unit vector in $\Bbb R^\ell$ and let $u$ be unit vector in 
$\Bbb R^{d-\ell}$ so that 
$$|
\inn {u}{\nabla_{x'}}\inn{v}{\nabla_{y'}}\vth\!\cdot\!S(0,0)|=\max_{U\in S^{d-\ell}}
|\inn {U}{\nabla_{x'}}\inn{v}{\nabla_{y'}}\vth\!\cdot\!S(0,0)|.
\tag 5.6
$$ Suppose further that the standard assumptions of \S2.1 and
(2.16) hold and 
$$\inn{u}{\nabla_{x'}}S(x,x')=0
\tag 5.7$$
for all $|x|\le \eps$. Let 
$a (\eta'')$ be supported in
 $\{\eta'':|\frac{\eta''}{|\eta''|}-\vth|\le \eps^5\}$ and satisfy
$|\partial^\alpha  a(\eta'')|\lc |\eta''|^{-|\alpha|}$ for 
all admissible multiindices $\alpha$.
Let $$\Theta= a(D).$$ 
Then the operator $\Theta \cB$  is bounded from
 $L^p$ to $L^q$ and its  operator norm satisfies the estimate
$$
\|\Theta \cB\|_{L^p\to L^q}
\lc 1+\|\cR\|^{1-\frac p 2}_{L^p\to L^q}.
$$
\endproclaim

\subheading{Proof of Proposition 5.3}

We can rewrite $\cB$ as 
$$\cB=\sum_{j\in J}\sum_{k>j}\omega(2^{-2k} |D''|)\cR_j.$$
Let
$\cL_k$ be defined by 
$$
\widehat {\cL_k  f} (\eta)= 
\omega(2^{-2k} |\eta''|)a(\eta'').
$$
then $\Theta\cB=\sum_{j\in J}\sum_{k>j} \cL_k \cR_j $.

We shall now introduce an angular Littlewood-Paley decomposition (as in 
\cite{14}) and 
proceed for the proof of our endpoint estimate using a well known argument by M. Christ
 (his preprint  \cite{3} is  unpublished but the argument has been used in 
various related articles on $L^p$ improving properties of convolution operators; 
for a rather general  formulation see \cite{10}).
Define operators
$P_{k,j}$, $\widetilde P_{k,j}$  by
$$
\align
&P_{k,j} = \sum_{i=-M}^M
\zeta(2^{-2k+j+i} |\inn{u}{D'}|)
\tag 5.8.1
\\
&\widetilde P_{k,j} =
\sum_{i=-M-10}^{M+10}
\zeta(2^{-2k+j+i} |\inn{u}{D'}|)
\tag 5.8.2
\endalign
$$
(we have chosen $2^M\ge c_0^{-1}\eps^{-10}$). Define also
$$
\widetilde L_k=
\sum_{i=-10}^{10} \omega(2^{-2k+i}|D''|) 
$$

The operator $\Theta\cB$  is then decomposed as
$$
\Theta\cB=
\sum_{j\in J}\sum_{k>j} \cL_k \cR_j  =\cT+  \cE_1+\cE_2+\cE_3
$$
where
$$
\align
\cT&=\sum_{j\in J} \sum_{k>j} \cL_k P_{k,j}
\cR_j \widetilde P_{k,j} \widetilde L_k
\tag 5.9
\\
\cE^1&=\sum_{j\in J} \sum_{k>j} \cL_k (I-P_{k,j})
\cR_j \widetilde P_{k,j} \widetilde L_k
\tag 5.10
\\
\cE^2&=\sum_{j\in J} \sum_{k>j} \cL_k 
\cR_j (I-\widetilde P_{k,j}) \widetilde L_k
\tag 5.11
\\
\cE^3&=\sum_{j\in J} \sum_{k>j} \cL_k 
\cR_j (I- \widetilde L_k).
\tag 5.12
\endalign
$$

The main term is represented by $\cT$, and we shall show that the 
operators $\cE^1$, $\cE^2$ and $\cE^3$  have quantitative properties 
similar to or better than the operator considered in Lemma   5.2.

For the  main term we use the known  argument in the translation invariant case \cite{3}.
Let $\cT_\vect$ denote the operator acting on $L^p(\ell^2(\bbZ^2))$ 
functions $F=\{F_{j,k}\}$ by
$$[\cT_\vect F]_{j,k}= \cR_j F_{j,k}.$$
By Littlewood-Paley theory and   complex interpolation (note that $p\le 2$)
$$\align
\|\cT \|_{L^p\to L^q} &\lc \|\cT_\vect\|_{L^p(\ell^2)\to L^q(\ell^2)}
\\
&\lc \|\cT_\vect\|_{L^p(\ell^p)\to L^q(\ell^p)}^{p/2}
 \|\cT_\vect\|_{L^p(\ell^\infty)\to L^q(\ell^\infty)}^{1-p/2}.
\tag 5.13
\endalign
$$
From  Proposition 3.2 and Minkowski's inequality it follows that
$$\|\cT_\vect\|_{L^p(\ell^p)\to L^q(\ell^p)}\lc 1.
\tag 5.14$$ 
Also by the pointwise  inequality 
$|\cR_j(f)|\le \cR(|f|)$ and the
positivity of  $\cR$ we have 
$$
\sup_{j,k}|\cR_j F_{j,k}(x)|\le
\cR[ \sup_{j,k} |F_{j,k}|](x)
$$ 
so that 
$$\|T_\vect\|_{L^p(\ell^\infty)\to L^q(\ell^\infty)}\lc 
\|\cR\|_{L^p\to L^q}.
\tag 5.15
$$

Therefore in view of Lemma 5.2 and (5.13-15)
$$
\|\cT\|_{L^p\to L^q}\le C(1+ 
\|\cR\|_{L^p\to L^q}^{1-p/2}+ \sum_{i=1}^3 \|\cE_i\|_{L^p\to L^q})
\tag 5.16
$$
Consequently the proof of Proposition 5.3 will be  complete once we verify
the uniform $L^p\to L^q$ boundedness of the operators  
$\cE^1$, $\cE^2$, $\cE^3$.

It will be convenient to work with 
oscillatory integral representations of the kernels of $\cR_j$. 
Since the Fourier transform of $c_{\ell,\r}|\cdot|^{\rho+\ell}$ is $|\xi|^\r$  (see \cite{8})
we can  write the kernel  $R_j$ of $\cR_j$ as an oscillatory integral
$$
R_j(x,y) =\kappa_j(x,y)|x'-y'|^{\s-d+\ell}\int
e^{\ic \inn {\tau} {y''-S(x,y')}} 
|\tau|^\r 
d\tau.
$$

For $k\ge 1$ we denote by $\cR_{j}^k$ the operator with integral kernel
$$
R^k_j(x,y)=\kappa_j(x,y)|x'-y'|^{-(d-\ell-\s)}
\int
e^{\ic \inn {\tau} {y''-S(x,y')}} 
\omega(2^{-2k} |\tau|)\zeta_0\big(\eps^{-4}|\frac{\tau}{|\tau|}-\vth|\big)
|\tau|^{\rho}d\tau;
$$
the operator 
$\cR_{j}^0$ is defined similarly but with 
$\omega(2^{-2k} |\tau|)$ replaced by $\omega_0(|\tau|)$.

\proclaim{Lemma 5.4} 
(i) The operator 
$\sum_j \cR^0_j$ maps $L^p$ to $L^q$.

(ii) Let $s\ge 0$. Let $Z_s(x,y)$ denote the distribution kernel of
the operator
$\sum_j \cL_{j+s}(\cR_j-\sum_{i=-4}^4\cR_{j}^{j+s+i})$. Then
$$
|Z_s(x,y)|\lc 4^{-s} |K_2^{\r,\s}(x,y)|
$$ 
where  
$K_2^{\r,\s}$ is defined in (4.3). Thus this operator maps $L^p\to L^q$ with 
operator norm $O(4^{-s})$.
\endproclaim

\demo{Proof}
(i) It is easy to see that  by the theorem on fractional integration the operator
$\sum_j \cR^0_j$ maps $L^p$ to $L^q$, provided that $1<p<q<\infty$ and $(d-\ell)(1/p-1/q)\le \s$.
However the condition $(d-\ell)(1/p-1/q)\le \s$ is implied by
$(d+\ell)(1/p-1/q)\le \s-2\r$ and $-\r(d-\ell)/\ell\le \s$ which is assumed throughout 
this section.

(ii)  Note that
$$
\cR_j-\sum_{i=-4}^4 \cR_j^{j+s+i}=
\sum_{r\ge 5} \cR^{j+s+r}_j+
\cS_{j,j+s}^0 +\sum_{r\ge -4}\cV_{j,j+s+r}+\cV_{j,j+s}^0
$$
where the kernels $S^0_{j,k}$, $V_{j,k}$ and $V^0_{j,k}$ 
of
$\cS^0_{j,k}$, $\cV_{j,k}$ and $\cV^0_{j,k}$ are given by
$$
\align
S^0_{j,k}(x,y)&=
\kappa_j(x,y)|x'-y'|^{-(d-\ell-\s)}
\int
e^{\ic \inn {\tau} {y''-S(x,y')}}
\omega_0(2^{-2(k-5)} |\tau|)
|\tau|^\rho
\zeta_0(\eps^{-4}|\frac{\tau}{|\tau|}-\vth|) 
d\tau
\\
V_{j,k}(x,y)&=
\kappa_j(x,y)|x'-y'|^{-(d-\ell-\s)}
\int
e^{\ic \inn {\tau} {y''-S(x,y')}}
\omega(2^{-2k} |\tau|)
|\tau|^\rho
\big(1- \zeta_0(\eps^{-4}|\frac{\tau}{|\tau|}-\vth|)\big)
d\tau
\\
V_{j,k}^0(x,y)&=
\kappa_j(x,y)|x'-y'|^{-(d-\ell-\s)}
\int
e^{\ic \inn {\tau} {y''-S(x,y')}}
\omega_0(2^{-2(k-5)} |\tau|)
|\tau|^\rho
\big(1- \zeta_0(\eps^{-4}|\frac{\tau}{|\tau|}-\vth|)\big)
d\tau.
\endalign
$$

We shall now show that 
the distribution kernel of
$\sum_j \cL_{j+s} \cR_{j}^{ j+s+r}$   is for $r\ge 5$  controlled by 
$4^{-(s+r)}
K_2^{\r,\s}$ ({\it cf.} (4.3)).
Also 
the  kernels  of
$\sum_j \cL_{j+s} \cS_{j, j+s}^0$ and 
$\sum_j \cL_{j+s} \cV_{j, j+s}^0$ are   
   bounded by 
$4^{-s}
K_2^{\r,\s}$; we shall omit the entirely analogous argument.

The kernel of $\cL_n \cR_{j}^k$ is given by
$$
\multline
K_{j,k,n}(x,y)=
(2\pi)^{-\ell}\,
\iiint e^{\ic [\inn{x''-w''}{\eta''}+\inn{\tau}{y''-S(x',w'',y')}]}
\,\times\,\\
\omega(2^{-2k}|\tau|) |\tau|^\rho
\omega(2^{-2n}|\eta''|)a(\eta'')
\zeta_0(\eps^{-5}|\frac{\tau}{|\tau|}-\vth|) \frac{\ka_j(x',w'',y',y'') }
{|x'-y'|^{d-\ell-\s}} dw'' \, d\eta''\, d\tau.
\endmultline
$$
We need to estimate this kernel when $k\ge n+5$, and $n\ge j$.
The $w''$-gradient of the phase function is 
$-\eta''-\nabla_{w''}(\tau\cdot S(w,y'))$
and since $\|S_{w''}-I_{\ell,\ell}\|\ll\eps^{1/2}$
this  gradient is now $\approx 2^{2k}$ (note that it would be  
$\approx 2^{2n}$ if we 
worked with
$\cL_n S_{j,n}^0$).

We use integration by parts with respect to $w''$ 
followed by integration by parts with respect to $\tau$ and $\eta$.
Observe that with each differentiation of
$\kappa_j(x',w'',y)$ we loose a factor of $2^{2j}$, 
the main contribution coming from
differentiating 
$\zeta_0(|w''-S(x',w'',y')|/|w'-y'|^2)$.
Thus we gain 
$2^{-2k+2j}$ with each integration by parts in $w''$.
As a result we obtain that the kernel of
$\cL_n \cR_{j}^k$ is dominated by a constant times 
$$
\align
& 2^{-(2k-2j)N_0}
\int |x'-y'|^{\s-d+\ell}\frac{2^{2k(\ell+\rho)}}{(1+2^{2k}|y''-S(x',w'',y')|)^{N_1}}
\frac{2^{2n\ell}}{(1+2^{2n}|x''-w''|)^{N_1}}  dw''
\\
&\lc \min\{2^{-(2n-2j)(N_0-N_1)}, 2^{-(2k-2j)(N_0-N_1)}\}
|x'-y'|^{\s-d+\ell}\frac{2^{2k(\ell+\rho)}}{(1+2^{2k}|y''-S(x',x'',y')|)^{N_1}};
\endalign
$$
here we choose $N_0\gg N_1$.
Moreover the kernel of the operator
$\cL_n \cS_{j,k}$ is of course supported where $|x'-y'|\approx 2^{-j}$. 
The asserted pointwise estimate for 
$\sum_j \cL_{j+s} \cR_{j}^{ j+s+r}$   is now a consequence of summing 
geometric series.

The same argument applies to the operators
$\sum_j \cL_{j+s} \cV_{j, j+s+r}$, $r\ge -4$.
Note that the above restriction  $r>4$ (or $k> n+4$) is not necessary  now 
in view of the factor
$(1- \zeta_0(\eps^{-4}|\frac{\tau}{|\tau|}-\vth|))$; namely  
the assumptions
$\eta''\in \supp a$ (hence $\big|\eta''/|\eta''|-\vth\big|\le \eps^{5} $) and
$|\tau/|\tau|-\vth|\ge \eps^{4}/2\gg\eps^{5}$
 guarantee that $|-\eta''-\tau\!\cdot\! S_{w''}(w,y')|\approx 
\max\{|\eta''|,|\tau|\}$ which is sufficient to carry out the 
above integration by parts arguments.
\qed
\enddemo

We shall now bound the operators $\cE^{1}$, $\cE^2$ and $\cE^3$  in (5.10-12).
However we first modify these operators by  replacing 
$\cL_k\cR_j$ in the definitions (5.10-12) by 
$\sum_{i=-4}^4\cL_k\cR_j^{k+i}$.
Let for $i=-4,\dots,4$
$$\align
\cE^1_{j,k,i}&= \cL_k (I-P_{k,j})
\cR_j^{k+i} \widetilde P_{k,j} \widetilde L_k
\tag 5.17.1
\\
\cE^2_{j,k,i}&= \cL_k 
\cR_j^{k+i} (I-\widetilde P_{k,j}) \widetilde L_k
\tag 5.17.2
\endalign
$$
and
$$
\cE^3_{j,k,i}= \cL_k 
\cR_j^{k+i} (I- \widetilde L_k)
\tag 5.18
$$
and let $$\widetilde \cE^{1,i}= \sum_{j\in J} \sum_{k>j} \cE^1_{j,k,i}, \quad i=-4,\dots,4;
\tag 5.19$$
 similarly define
$\widetilde \cE^{2,i}$, $\widetilde \cE^{3,i}$.

\proclaim{Lemma 5.5}
The operators 
 $\cE^1-\sum_{i=-4}^4
\widetilde \cE^{1,i}$,
 $\cE^2-\sum_{i=-4}^4\widetilde \cE^{2,i}$, and
 $\cE^3-\sum_{i=-4}^4\widetilde \cE^{3,i}$ are bounded from $L^p$ to $L^q$.
\endproclaim
\demo{ Proof}
This is a  consequence of Lemma 5.4. We use it in conjunction with
 Littlewood-Paley theory, the iterated version of the 
Fefferman-Stein vector-valued maximal function  and
the Marcinkiewicz-Zygmund theorem on vector-valued extensions  
of $L^p\to L^q$ bounded operators (\cite{7}, \cite{22}).
We use the pointwise estimate $|P_{k,j}g|\le \fM g$
where $\fM$ denotes the strong maximal function.
Let $F^{\r,\s}$ be the fractional integral operator with distribution kernel
$K_2^{\r,\s}$.
Then
$$
\align
\|\cE^1 f&-\sum_{i=-4}^4
\widetilde \cE^{1,i}f\|_q\lc
\sum_{s\ge 0}\Big\|\Big(\sum_{j\in J}|(I-P_{j+s,j})
\cL_{j+s}
\big(\cR_j-\sum_{i=-4}^4 \cR_j^{j+s+i}\big)\widetilde P_{j+s,j}\widetilde L_{j+s} f|^2\Big)^{1/2}\Big\|_q
\\
&\lc
\sum_{s\ge 0}4^{-s}\Big\|\Big(\sum_{j\in J}\big[\fM
F^{\r,\s}[|\widetilde P_{j+s,j}
\widetilde L_{j+s} f|]\big]^2\Big)^{1/2}\Big\|_q
\lc
\sum_{s\ge 0}4^{-s}\Big\|\Big(\sum_{j\in J}\big[
F^{\r,\s}[|\widetilde P_{j+s,j}\widetilde L_{j+s} f|]\big]^2\Big)^{1/2}\Big\|_q
\\
&\lc
\sum_{s\ge 0}4^{-s}\Big\|\Big(\sum_{j\in J}|
\widetilde P_{j+s,j}\widetilde L_{j+s} f|^2\Big)^{1/2}\Big\|_p \lc \|f\|_p.
\endalign
$$
The other estimates are proved in a similar way.\qed
\enddemo

As a consequence of Lemma 5.5 it remains,  in order to conclude the proof of Proposition 5.3,
to show that   the operators
 $\widetilde \cE^{1,i}$,
 $\widetilde \cE^{2,i}$,
 $\widetilde \cE^{3,i}$ are bounded from $L^p$ to $L^q$.
We shall show that
 $\widetilde \cE^{1,i}$ maps $L^p$ to $L^q$. The proof of the boundedness 
of $\widetilde \cE^{2,i}$ is very similar and will therefore be  omitted.
Finally, the arguments in the proof of Lemma 5.4 show  the 
$L^p\to L^q$ boundedness  of $\widetilde  \cE^{3,i}$; the details will be  omitted as well.

\demo{Boundedness of  $\widetilde \cE^{1,i}$}
We   analyze the kernel  of $\cL_k(I-P_{k,j})\cR_j^{k+i}$
which is given by
$$
K_{k,j,i}(x,y)= (2\pi)^{-\ell-1}\iint\iiint
e^{\ic \varphi(x,t,h'',y,\tau,\lambda,\eta'')}
a_{k,j,i}(x,t,h'',y,\tau,\lambda,\eta'') \,d\tau d\eta'' d\lambda \,dh'' dt
\tag 5.20
$$
where
$$
\varphi(x,t,h'',y,\tau,\lambda,\eta'')
=-t\lambda-\inn{\eta''}{h''}-\inn{\tau}{S(x'+tu,x''+h'',y')-y'' }
\tag 5.21
$$
and
$$\multline
a_{k,j,i}(x,t,h'',y,\tau,\lambda,\eta'')=
a(\eta'')\omega(2^{-2k}|\eta''|)\omega(2^{-2(k+i)}|\tau|) |\tau|^{\rho}
\\
\zeta_0(\eps^{-1}_2|\tau/|\tau|-\vth|)
\chi(x'+tu,x''+h'',y)\kappa_j(x'+tu,x''+h'',y)|x'+tu-y'|^{\s-d+\ell} 
(1-\zeta_M(2^{-2k+j}|\lambda|))
\endmultline
\tag 5.22
$$
with $\zeta_M=\sum_{s=-M}^M\zeta(2^{s}\cdot)$.

\demo{ Claim}  {\it For $s\ge 0$, $i=-4,\dots,4$  we have
$$
|K_{j+s,j,i}(x,y)|\lc 4^{-s} |K_2^{\r,\s}(x,y)|
$$
uniformly in $j$.
Here the right  hand side  is defined in (4.3).}  

\enddemo

Taking  the claim for granted we can argue 
as in the proof of  Lemma  5.5 and obtain  using Littlewood-Paley theory and the boundedness of the operator $F^{\r,\s}$ with kernel $K_2^{\r,\s}$
$$
\align
\big\|\widetilde E^{1,i} f\big\|_q&=\Big\|\sum_{s>0}\sum_{j\in J}
\Cal L_{j+s}(I-P_{j+s,j}) \cR_j^{j+s+i} \widetilde P_{j+s,j} \widetilde L_{j+s} f
\Big\|_q
\\
&\lc \sum_{s>0}\Big\|\Big(\sum_{j\in J}\big|
\Cal L_{j+s}(I-P_{j+s,j}) \cR_j^{j+s+i} \widetilde P_{j+s,j}
 \widetilde L_{j+s} f
\big|^2\Big)^{1/2}\Big\|_q
\\
&\lc \sum_{s>0}4^{-s}\Big\|\Big(\sum_{j\in J}\big|
F^{\r,\s} [|\widetilde P_{j+s,j} \widetilde L_{j+s} f|]
\big|^2\Big)^{1/2}
\Big\|_q
\\
&\lc \sum_{s>0}4^{-s}\Big\|\Big(\sum_{j\in J}\big|
\widetilde P_{j+s,j} \widetilde L_{j+s} f
\big|^2\Big)^{1/2}
\Big\|_p \lc \|f\|_p
\endalign
$$

We proceed to prove the pointwise estimate claimed above.
We note that 
$$
a_{k,j,i}(x,t,h'',y,\tau,\lambda,\eta'')=0 \quad\text{ if } |\lambda|\in
[2^{2k-j-M+4},2^{2k-j+M-4}].
\tag 5.23
$$
Now we  first  integrate by parts many times 
 in (5.20) with respect to $t$;
this is 
then followed by an integration by parts in the $(\lambda, \eta'', \tau)$ variables.

Note that because of $\inn{u}{\nabla_{w'}  S(y',w'',y')} =0$  we may expand
$$
\align
\partial_t
& \varphi(x,t,h'',y',\tau,\lambda,\eta) =
-\lambda-\inn{u}{ \tau\!\cdot\! S_{x'}(x'+tu,x''+h'',y')} 
\\
&=
-\la+\inn{u}{\tau\!\cdot\! S_{x'x'}(y',x''+h'',y')(x'+tu-y')} + 
\tau\!\cdot\! r_1(x,y',t,h'')
\\
&=-\lambda+\inn{u}{\tau\!\cdot\! S_{x'x'}(0,0,0)(x'+tu-y')} + \tau\!\cdot
\! \big( \sum_{\nu=1,2}r_\nu(x,y',t,h'')\big)
\tag 5.24
\endalign
$$ where 
$$
\align
&|r_1(x,y',t,h'')|\le C_1 |y'-x'-tu|^2 
\\
&|r_2(x,y',t,h'')|\le C_1\eps^{10} |y'-x'-tu|.
\endalign 
$$
Differentiating (5.7) we see that 
$$\inn{u}{\widetilde S_{x'x'}(x,x')
+ \widetilde S_{x'y'}(x,x')}
=0
$$
and by (2.7-8) 
and the  choice of $u$  we deduce that
$$c_0 2^{2k-j-2}\le \big|\inn{u}{\tau\!\cdot\! S_{x'x'}(0,0,0)(x'+tu-y')}\big|
 \le c_0^{-1} 2^{2k-j+3}
$$
and consequently,
by our choice of $M$
$$2^{2k-j-M+5}\le c_02^{2k-j-2}\le \big|
\partial_t \varphi(x,t,h'',y',\tau,\lambda,\eta'') +\lambda\big|
 \le c_0^{-1} 2^{2k-j+3}
\le 2^{2k-j+M-5}
$$
on the support of the symbol; hence by (5.23)
$$|\partial_t \varphi(x,t,h'',y',\tau,\lambda,\eta'') | \gc \max\{\lambda,\,
 2^{2k-j}\}.
$$
Moreover the higher derivatives of the phase 
functions are $O(2^{2k-j})$. 
Taking $s$ derivatives of  $\ka_j$  with respect to $w'$ (in any direction)
 causes a blowup of size $O(2^{2js})$ which 
would  be  too much for our argument.
Fortunately, in view of the assumption $\inn{u}{\nabla_{x'}  S(y',w'',y')} =0$ 
we have the better estimate
$$
(\inn{u}{\nabla_{w'})^s\kappa_j(w,y))} =O(2^{js}).
$$

Thus we may perform integration by parts in the $t$ variables 
and gain factors  of size $2^{(2j-2k)N}$. This is then  followed by an integration by parts 
in the frequency  variables and we obtain 
$$
\multline
|K_{k,j,i}(x,y)|\lc
2^{j(d-\ell-\sigma)} 2^{2k\r}
2^{-(2k-2j) N_1}
 \int \chi_j(x'+tu-y')
\chi(x'+tu,x''+h'',y) \,\times
\\
\frac{2^{2k-j }}{(1+2^{2k-j}|t|)^{N_2}}
\frac{2^{2k\ell}}{(1+2^{2k}|h''|)^{N_2}}
\frac{2^{2k\ell}}{(1+2^{2k}|y''-S(x'+tu,x''+h'',y')|)^{N_3}} dt dh''.
\endmultline
$$
Now observe that 
$$
|S(x'+tu,x''+h'',y')-S(x,y')|\lc |h''|+2^{-j}|t|+|t|^2
$$ and therefore 
$$
\frac{2^{2k\ell}}{(1+2^{2k}|y''-S(x'+tu,x''+h'',y')|)^{N_3}} \lc
\frac{2^{2k\ell}} {(1+2^{2k}|y''-S(x,y')|)^{N_3}} 
\big(1+ 2^{2k-j}|t|+|t|^2
+2^{2k}|h''|)^{N_3}
$$
This yields 
$$
\multline
|K_{k,j,i}(x,y)|\lc 
2^{-(2k-2j) (N_1-\r-\ell)}
 2^{j(d+\ell-\sigma+2\rho)}  
(1+2^{2k}|y''-S(x,y')|)^{-N_3}\,\times\, 
\\
\iint_{\bbR\times\bbR^\ell} \chi_j(x'-y'-tu)
\chi(x''+tu,x''+h'',y)
\frac{2^{2k-j } }{(1+2^{2k-j}|t|)^{N_2-N_3}}
\frac{2^{2k\ell}  
}{(1+2^{2k}|h''|))^{N_2-N_3}}
dt  dh''
\endmultline$$
where $\chi_j$ denotes the characteristic function of 
$[2^{-j-1},2^{-j+1}]
\cup[-2^{-j+1},-2^{-j-1}]$.

This integral is straightforward to estimate.
Observe that
$ 2^{j(d+\ell-\sigma+2\rho)}  
(1+2^{2k}|y''-S(x,y')|)^{-N_3}$ is bounded by $ |y''-S(x,y')|^{-(d+\ell-\sigma+2\rho)/2}  $;
thus if 
$|x'-y'|\le C2^{-j}$ we use either this bound or 
the bound  $ 2^{j(d+\ell-\sigma+2\rho)} $ and estimate
$|K_{k,j,i}(x,y)|$ by $C 
2^{-(2k-2j) (N_1-\r-\ell)} K_2^{\r,\s}(x,y)$.

Next, if $C2^{-j}\le|x'-y'|\le \eps$ and  $|y''-S(x,y)|\le \eps$
 then
$\chi_j(x'-y'-tu)$ vanishes unless $|t|\ge c |x'-y'|$. In 
this case the contribution of the   $t$ integral above  is  
$$O\big((2^{j-2k}|x'-y'|^{-1})^{N_2-N_3-1}\big)  
+
O\big((2^{-2k}|x'-y'|^{-1})^{N_2-N_3-d+\ell}\big).$$

Thus in this case
$$|K_{k,j,i}(x,y)|\lc 
2^{-(2k-2j) (N_1-\r-\ell)} 2^{j(d+\ell-\sigma+2\r)}    (1+ 2^j|x'-y'|)^{-2N} 
(1+2^{2k}|y''-S(x,y')|)^{-N} 
$$
where
$2N= \min \{N_2-N_3-d+\ell, N_3\}$. We may choose 
$2N\le N_1+2d$ and $ N\gg d$ and again  the bound 
$|K_{k,j,i}(x,y)|$ by $C 
2^{-(2k-2j) (N_1-\r-\ell)} K_2^{\r,\s}(x,y)$ is straightforward.
Thus we have established the pointwise estimate claimed above. 
This concludes the proof of Proposition 5.3.
\enddemo

\subheading{Proof of Theorem 5.1, conclusion}
We have to prove that  $\cR$ in (5.4) maps $L^p$ to $L^q$; 
assuming the angular  localization (5.1)  in the 
$x'-y'$ variables.
We split the identity operator  as 
$E_0+\sum_\nu \Theta_\nu$  where $E_0=\eta_0(D'')$ and $\eta_0$ is compactly supported 
in $\{\eta'':|\eta''|\le 1000\}$. Moreover let $\Theta_\nu=a_\nu(D'')$ where $a_\nu$ is a 
constant coefficient symbol of order $0$ supported in  
$$\{\eta'': |\frac{\eta''}{|\eta''|}-\vartheta_\nu|\le \eps^5, |\eta''|\ge 100
\};$$
we can arrange this decomposition so that the sum in $\nu$ is extended over
 $O(\eps^{-5(\ell-1)})$ terms. Clearly it suffices to bound $E_0 \cR$ and 
$\Theta_\nu\cR$
for all $\nu$.
We first note that the argument of Lemma 5.2 shows that
$E_0 \cR$  maps $L^p\to L^q$ if $(d+\ell)(1/p-1/q)\le \sigma-2\r$.

It remains to consider $\Theta_\nu\cR^\sigma$ for fixed $\nu$.
Let $u_\nu$ be a unit vector in $\Bbb R^\ell$ so that
$$|
\inn {u_\nu}{\nabla_{x'}}\inn{v}{\nabla_{y'}}\vth_\nu\!\cdot\!S(0,0)|=\max_{U\in S^{d-\ell}}
|\inn {U}{\nabla_{x'}}\inn{v}{\nabla_{y'}}\vth_\nu\!\cdot\!S(0,0)|.
$$ 
Now denote by  $\cQ^\nu$ the change of variable $\cQ(\cdot,u_\nu)$ as defined in
\S2.2, moreover define $\fQ_\nu h(w)= h(\cQ_\nu w)$ for functions supported in $B_{\eps^9}$.
Let $\cR^{\nu} =\fQ_\nu \cR \fQ_\nu^{-1}$; then the assumptions of Proposition 5.3 apply to 
$\cR^\nu$ (with $u=u_\nu$).

Define $\widetilde \Theta_\nu=\widetilde a_\nu(D'')$ so that
$\widetilde a_\nu$ is  supported in  
$\{\eta'': |\frac{\eta''}{|\eta''|}-\vartheta_\nu|\le \eps^2, |\eta''|\ge 10\};$
and $\widetilde a_\nu(\eta'')=1$
if $|\frac{\eta''}{|\eta''|}-\vartheta_\nu|\le \eps^2$ and $|\eta''|\ge 20\}$.
Then by Proposition 5.3 and Lemma 5.2

$$
\|\widetilde \Theta_\nu \cR^\nu \|_{L^p\to L^q}
\le C(1+\|\cR^\nu\|^{1-\frac p2}_{L^p\to L^q})
\tag 5.26
$$
But in view of the support properties of the kernel of $\cR$ and the local $L^p$ and $L^q$ boundedness of the operators $\fQ_\nu$ and 
$\fQ_\nu^{-1}$  we get 
$$
\|\cR^\nu\|_{L^p\to L^q} \lc
\|\cR\|_{L^p\to L^q}.
$$

To conclude the proof we split
$$
\align
\cR &
=E_0 \cR+\sum_{\nu}\Theta_\nu \fQ_\nu^{-1}\cR^\nu\fQ_\nu
\\&=E_0 \cR+\sum_{\nu}\Theta_\nu \fQ_\nu^{-1}
\widetilde \Theta_\nu \cR^{\nu}\fQ_\nu
+\sum_{\nu}\Theta_\nu \fQ_\nu^{-1} (I-
\widetilde \Theta_\nu) \cR^{\nu}\fQ_\nu.
\endalign
$$
By  (5.26)
$$
\|\Theta_\nu \fQ_\nu^{-1}
\widetilde \Theta_\nu \cR^{\nu}\fQ_\nu\|_{L^p\to L^q}\lc 1+\|\cR\|_{L^p\to L^q}^{1-p/2}
\tag 5.27
$$
and it remains to show that
$$\|\Theta_\nu \fQ_\nu^{-1} (I-\widetilde \Theta_\nu) \cR^{\nu}\fQ_\nu\|_{L^p\to L^q}\lc 1.
\tag 5.28
$$

Now let $L_0=\omega_0(|D''|)$ and $L_k=\omega(4^{-k}|D''|)$. 
We analyze the kernel of $L_k\Theta_\nu \fQ_\nu^{-1} (I-\widetilde \Theta_\nu)L_{k'}$, denoted by
$H_{k,k',\nu}(x', x'',y'')$.
The inverse change of variable $\cQ_\nu^{-1}$ 
is of the form $x\mapsto (x', \cG_\nu(x))$, with $\|(\cG_{\nu})_{x''}-I_{\ell,\ell}\|\le
\eps^7$ ({\it cf.} (2.17/18)).
Thus 
$H_{k,k',\nu}$ 
is given by
$$
\multline H_{k,k',\nu}(x',x'',y'')=\\
\iiint
e^{\ic(\inn{x''-z''}{\eta''}+\inn{\cG_\nu(x',z'')-y''}{\xi''})} \omega(4^{-k}|\eta''|)
\omega(4^{-k'}|\xi''|)a_\nu(\eta'') (1-\widetilde a_\nu(\xi''))
dz'' d\xi'' d\eta''
\endmultline
$$
The $z''$-gradient of the phase function is of size $\approx \max\{4^k, 4^{k'}\}$, therefore we may argue as in the proof of Lemma 5.4 above.
In particular, after additional integration by parts in $\xi'', \eta''$ when $x$  is large we obtain that
$$|H_{k,k',\nu}(x',x'',y'')|
\lc \min\{ 4^{-kN_1}, 4^{-k' N_1}\}(1+|x''|)^{-N_2}.$$

In view of the localization properties of $\cR^\nu$ and the $L^p$ 
boundedness of $\cR^\nu$ it follows that
$$\|L_k\Theta_\nu \fQ_\nu^{-1} (I-\widetilde \Theta_\nu)L_{k'} \cR^{\nu}\fQ_\nu\|_{L^p\to L^q}
\lc \min\{4^{-k}, 4^{-k'}\}
$$
and as a consequence (5.28) holds.

Putting all the estimates together we obtain that
$$
\| \cR\|_{L^p\to L^q}
\lc 1+\|\cR\|_{L^p\to L^q}^{1-p/2}
\tag 5.29
$$
and since we already know the finiteness of 
$
\| \cR\|_{L^p\to L^q}$ the estimate (5.29) implies a bound uniform in the family $J$.\qed

We can now give the 

\demo{\bf Proof of Theorem 1.2.3}
By summing geometrical series we see from Lemma 4.3 that the operator 
$\sum_{m\ge 0}\sum_{k>m} T_{k,m}$
 can be pointwise bounded by a combination of operators handled in
Theorem 5.1; in this calculation we use that $\rho$ is negative. Moreover 
the operator 
$\sum_{m\ge 0}\sum_{k\le m} T_{k,m}$ is bounded by Proposition 4.4. The assertion 1.2.3  follows.
\qed
\enddemo

\subheading{Necessary conditions}
The necessity of the conditions in Theorems 1.1  and 5.1 follows from standard examples.
For the sake of completeness we shall briefly describe them.
We assume that  $\rho\le 0$ and $1\le p\le q\le \infty$ and consider the operator $R^{\r,\s}$.
We remark that
for the case $\rho>0$, the conditions in 1.2.1 also cannot be improved. This is because any
 strict 
improvement would yield to an improvement in the case $\rho=0$, by interpolation with the 
estimates for a negative $\r_1$ close to $0$. 

Let $B_\delta$ be the ball of radius $\delta\ll\eps^{10}$, centered at the origin, and let 
$\chi_\delta$ be the characteristic function of $B_\delta$.
Then $\|\chi_\delta\|_p\gc \delta^{d/p}$ and $R^{\r,\s}\chi_\delta\gc \delta^{d-\ell-\r}$
 on the set $\{x:|x'|\le \delta^2, |x''-\fS(0,x')|\le c\delta\}$ for small $c$.
Thus 
$\|R^{\r,\s}\chi_\delta\|_q\gc \delta^{d-\ell-\r-(d-\ell)/q}$ and we see that the 
condition $d/p-\ell/q\le d-\ell-\r$ is necessary.
By applying the same example to the adjoint operator we get the necessary condition
$\ell/p-d/q\le -\r$. Thus 
$(1/p, 1/q)$ belongs to the pentagon with  corners
$(1,1)$, $(0,0)$, $(1,\frac{\rho+\ell}{\ell})$, $(\frac{-\rho}{\ell},0)$ and
$(\frac{d-\r}{d+\ell}, \frac{\r+\ell}{d+\ell})$ and this pentagon becomes the triangle
in Theorem 1.1 when $\r=0$.

If $\rho<0$ then the operator $R^{\r,\s}$ is not bounded from $L^1$ 
to $L^{{\ell}/{\ell+\rho}}$
as one checks that 
   one has the 
lower bound
$R^{\r,\s}\chi_\delta\gc \delta^d |x''-\fS(0,x')|^{-\rho-\ell}$
if
$C\delta\le |x''-\fS(0,x)|\le \eps $, with $C$ large.  By applying this to the 
adjoint operator it follwos that $R^{\r,\s}$ is not 
bounded from $L^{-\rho/\ell}$ to $L^\infty$.

Next let $P_\delta$ be the plate $\{y:|y'|\le \delta, |y''|\le \delta\}$ and let $f_\delta$ be the characteristic function of 
$P_\delta$, thus $\|f_\delta\|_p\lc \delta^{(d+\ell)/p}$.
One checks that in a  fixed fraction of $P_\delta$ one has the lower bound 
$R^{\r,\s}f_\delta(x)\gc \delta^{\s-2\r}$; in this calculation we use (2.2) and (2.6). Thus
$\|
R^{\r,\s}f_\delta\|_q\gc \delta^{\s-2\r+(d+\ell)/q}$ and  the condition 
$(d+\ell)(1/p-1/q)\le \s-2\r$ is necessary.
This concludes the proof  of necessity in Theorem 1.1 and Theorem 5.1.

A third necessary condition for the $L^p\to L^q$  boundedness of $R^{\r,\s}$ is 
$(d-\ell)(1/p-1/q)\le \s$. To see this
let $g_\delta$ be the  
characteristic function  of $\{y: |y'|\le \delta, |y''|\le \eps\}$.
Then 
$R^{\r,\s} g_\delta\ge \delta^\s$ for all $x$ in a fixed  fraction of this set
and from this one deduces the necessity of the condition
$(d-\ell)(1/p-1/q)\le \s$.
Notice that the condition
$(d-\ell)(1/p-1/q)\le \s$ is more restrictive than
$(d+\ell)(1/p-1/q)\le \s-2\r$  if and only if $\s<-\r(d-\ell)/\ell$; thus 
this example is only relevant to show the sharpness of 1.2.4.

\head{\bf 6. $\boldkey L^{\boldkey p}$ estimates 
for Fourier integral operators}
\endhead

It will be convenient to introduce some normalized classes of symbols.

Let $k>0$ and $0<m<k$. Then we denote by $\cS_{k,m}$ the class of symbols
$a(x,y,\xi,\tau)$  supported in
$$\alignedat2 &\{ (x,y,\tau,\xi): |x|+|y|\le \eps, 
2^{2k-1}\le|\tau|\le 2^{2k+1},
2^{m-1}\le|\xi|\le 2^{m+1}\} \quad &&\text{if } 0<m<k,
\\
&\{ (x,y,\tau,\xi): |x|+|y|\le \eps, 
2^{2k-1}\le|\tau|\le 2^{2k+1}, |\xi|\le 2\} \quad &&\text{if } m=0.
\endalignedat
\tag 6.1
$$
for which  (1.11) holds, with $\rho=\s=0$. Moreover, if $m>0$ let 
$\Sigma_m$ be the class of symbols 
$a(x,y,\xi,\tau)$  supported in
$$\aligned &\{ (x,y,\tau,\xi): |x|+|y|\le \eps, 
|\tau|\le 2^{2m+1},
2^{m-1}\le|\xi|\le 2^{m+1}\}
\endaligned
\tag 6.2
$$
such that (1.11) holds with $\r=\s=0$.

We recall that $\cT[a]$ denotes the  integral operator with kernel  (1.10).

\subheading{$\boldkey L^{\boldkey 2}$ estimates}

We shall assume that $a\in \cI^{\rho,-\sigma}$ and begin by proving 
$L^2$ estimates. These are quick consequences of what is already proved in
\cite{11}, and we shall be brief. It is shown in 
 in \cite{11} that 
$L^2$ boundedness holds if $2\rho -\sigma\le 0$, 
$0\le \sigma<d-\ell$. While the endpoint estimate corresponding to
$(\rho,\sigma)=((d-\ell)/2, d-\ell)$ may fail 
the proof of the estimates in \cite{11} still  provides  useful information 
which will be used in an interpolation argument in \S7.

\proclaim{Lemma 6.1} (i) 
Let   $a_m\in\Sigma_{m} $ and suppose that 
$\sup_{m\ge 1}|c_m|\le 1$.
Then $\sum_{m=1}^\infty c_m \cT[a_m]$ is bounded on $L^2$.

(ii) Let 
  $a\in S^{\frac{d-\ell}2,\ell-d}$ and suppose that $a(x,y,\tau,\xi)=0$ if 
$|\tau|\ge |\xi|^2$. Then $\cT[a]$ is bounded on $L^2$.
\endproclaim

\demo{\bf Proof} We note that the 
phase function 
$\Phi(x,y,\xi,\tau)=\inn{\xi}{x'-y'}+\inn{\tau}{y''-S(x,y')}$ parametrizes 
the diagonal  in $T^*\Omega\times T^*\Omega$ as a Lagrangian manifold; that is
$\{(x,\Phi_x,y,-\Phi_y):\Phi_\xi=0, \Phi_\tau=0\}$ is a subset of $\{(x,\xi,x\xi)\}$.

Because 
of the support restriction of $a_m$  the symbol $\sum_{m> 0}
c_m a_m $ belongs to the Calder\'on-Vaillancourt  symbol class
$S^{0}_{1/2,1/2}$.
It is shown in
the proof of Proposition 2.7
in  \cite{11} that H\"ormander's
 equivalence of phase function theorem remains valid with 
$S^{0}_{1/2,1/2}$ symbols and that consequently 
$\sum_{m=1}^\infty c_m \cT[a_m]$ is  a 
pseudodifferential operator of order $0$,
with symbols of type 
$(1/2,1/2)$. Thus the  $L^2$ boundedness follows from the 
Calder\'on-Vaillancourt theorem.
(ii) is an immediate consequence of (i).\qed
\enddemo


\proclaim{Lemma 6.2}
(i)
Let $m_0\ge 0$ be fixed and for $k> m_0$ let $m(k)$ be an integer such that $m_0\le m(k)<k$. 
Suppose that $\sup_{k\ge 1}|c_k|\le 1$ and that   
$\alpha_{k}\in
\cS_{k,m(k)}$.

Then the operator 
$\sum_{k>m_0} c_k 4^{k\frac{d-\ell}{2}} 
2^{-m(k)(d-\ell)} \cT[\alpha_{k}]$
 is bounded on $L^2 $, with norm independent of the  chosen sequence $\{m(k)\}$.

(ii) Suppose  $a\in S^{\frac{d-\ell}2,\ell-d}$ and suppose that
$a(x,y,\tau,\xi)=0 $ if $|\tau|\le C |\xi|^{1/2}$, and, for $m>0$,  let
 $a_m(x,y,\tau,\xi)=
\zeta(2^{-m}|\xi|) a(x,y,\tau,\xi)$. Then $\cT[a_m]$ is bounded on $L^2$ with 
operator norm independently  on $m$.

(iii) Let $\{\alpha_{k}\}$ be as in (i) and let 
$\eta\in S_{1/2,1/2}^0(\Omega\times\Omega, \Bbb R^d)$. Then 
the statement in (ii) remains valid if $\alpha_{k}$ is replaced by 
$\eta\alpha_{k}$.
\endproclaim

\demo{\bf Proof} For (i) we note that the 
kernel of  $\cT[\alpha_k]$ is given by 
$$
\int e^{\ic\inn{\tau}{y''-S(x,y')}} b_{k}(x,y,\tau) \, d\tau
\tag 6.3
$$
where
$$b_{k}(x,y,\tau)=\int \alpha_k(x,y,\tau,\xi) 
e^{\imath\inn{x'-y'}{\xi}}\, d\xi.\tag 6.4$$
Note that for every $k$ the $\xi$ integration is 
extended over a dyadic annulus $\{\xi:|\xi|
\approx 2^{m(k)}\}$ and thus 
$|b_{k,m}(x,y,\tau)|\lc 4^{k(d-\ell)/2}\approx |\tau|^{(d-\ell)/2}$.
Moreover, by examining the derivatives of $b_{k,m}$  one checks as in \cite{11} that  
$b_{k}$ is a symbol of order $(d-\ell)/2$ and type $(1/2,1/2)$.
Since the phase function involves $\ell$ frequency variables one may argue as in\cite{11}
and deduce that
$\sum_{k\ge m_0}c_k \cT[\alpha_{k}]$ are 
Fourier integral operators of order 
$0$ and type $(1/2,1/2)$,  hence bounded in $L^2$ (with bounds independent 
of the 
sequence $\{\alpha_k\}$).

Part (ii)  follows from part (i) with the choice $m(k)=m$  if we observe that 
the symbols
$a_m$ with the assumed support property can be decomposed as 
$C\sum_{k>m} 2^{k(d-\ell)/2} 2^{m(\ell-d)} c_{k,m} a_{k,m}$ where
$c_{k,m}\le 1$ and $a_{k,m}\in \cS_{k,m}$. Clearly the above  argument also 
proves (iii).\qed
\enddemo

\remark{Remark} The variant (iii) is included  in order to cover 
localizations of the form
$a_{k,m}(x,y,\tau,\xi)\zeta(2^j(|x'-y'|))$ if $j\le k$; these are of 
type $(1/2,1/2)$ since
$a_{k,m}$ is supported where $\tau\approx 2^{2k}$.
\endremark

\subheading{$\boldkey H^{\boldkey 1}\boldsymbol \to \boldkey L^{\boldkey 1}$
 estimates}

\proclaim{Lemma 6.3}
 Suppose  $0\le \s<d-\ell$,   $a\in S^{0,-\s}$,
 and suppose that
$a(x,y,\tau,\xi)$ is supported where $|\xi|\ge \frac 12|\tau|^{1/2}$. Let 
$$a_m(x,y,\tau,\xi)=\cases
a(x,y,\tau,\xi)\zeta(2^{-m}|\xi|) \quad &\text{if } m>0
\\
a(x,y,\tau,\xi)\zeta_0(|\xi|) \quad &\text{if }m=0
\endcases . \tag 6.5
$$
Then $\cT[a_m]$ maps $L^1$ boundedly to $L^1$, with operator norm $O((1+m)2^{-m\s})$.
\endproclaim

\demo{\bf Proof} 
The kernel $K_m$ can be written as $\sum_{k\le m} K_{k,m}$ where
$K_{k,m}$ is as in (4.6) and  satisfies (4.7) with $\rho=0$. The operator with kernel $K_{k,m}$ is clearly bounded on $L^1$, with norm $O(2^{-m\s})$.
\qed\enddemo

\proclaim{Lemma 6.4}
 Suppose  $a\in S^{0,-\s}$, $0\le \s<d-\ell$ and suppose that
$a(x,y,\tau,\xi)$ is supported where $|\xi|\le 2|\tau|^{1/2}$. 
 Let $a_m$ be as in (6.5).
Then $\cT[a_m]$ maps $H^1$ boundedly to $L^1$, with operator norm dominated by
$C2^{-m\s}$.
\endproclaim

\demo{Proof}
By the theorem on the atomic decomposition (\cite{7}, \cite{22})
it suffices to 
estimate $\cT[a_m]f_Q$ where
$f_Q$ is an $L^2$ function supported on a cube $Q$ with center $y_Q$ and 
sidelength $\delta_Q\ll 1$ so that 
$\|f_Q\|_2\le \delta_Q^{-d/2}$ and   $\int f_Q dx=0$.

We define the exceptional set 
$$W_Q=\{x:|x'-(y_Q)'|\le \eps, |x''-\fS(y_Q,x')|\le C\delta_Q\};$$
for large but fixed $C$; on this set   we shall use a  mixed norm $L^{1}(L^2)$ estimate.


We define
phase functions and amplitudes on $\Bbb R^\ell$ depending on the
 parameters $x',y'$.
Let
$$
b^{x',y'}_m(x'',y'',\tau)=\int a_m(x',x'',y',y'', \tau,\xi) 
e^{\ic\inn{x'-y'}{\xi}} d\xi
$$
and
$$
\Phi^{x',y'}(x'',y'',\tau)=\inn{\tau}{S(x',x'',y')-y''}.
$$
Denote by $\cT^{x',y'}_m$ the operator with kernel
$$K_m^{x',y'}(x'',y'')
=\int e^{\ic \Phi^{x',y'}(x'',y'',\tau)} b^{x',y'}_m(x'',y'',\tau)  d\tau.
$$
By an integration by parts one 
sees that 
$$
|\partial_{x'',y''}^\alpha \partial_\tau^\beta b_m^{x',y'}|
\le C_{\alpha,\beta} 
\frac{2^{m(d-\ell-\s)}}{(1+2^{m}|x'-y'|)^N}
$$
and  by the standard theory for pseudodifferential operators and their 
behavior under changes of variables it follows that
$$\big\|\cT_m^{x',y'}\big\|_{L^2(\bbR^\ell)\to L^2(\bbR^\ell)}
\lc 
\frac{2^{m(d-\ell-\s)}}{(1+2^{m}|x'-y'|)^N}.
$$

We now estimate the contribution on $W_Q$.
For fixed $x'$ set
$W_Q^{x'}=\{x'': (x',x'')\in W_Q\}$.
Let $f^{y'}(y'')=f(y',y'')$, then
$$\cT_m f_Q(x',x'')=\int_{y'}\cT_m^{x',y'}f^{y'} dy'.$$
On $W_Q$ we bound
$$\align
\int_{W_Q}|\cT_m  f_Q(x)| dx&\le
\int\limits\Sb|x'-(y_Q)'|\\ \le \eps\endSb
\int\limits\Sb |x''-\fS(y_Q,x')|\\ \le C\delta_Q\endSb
\int
|\cT_m^{x',y'}f_Q^{y'}(x')| dy' dx'' dx'
\\
&\lc\delta_Q^{\frac{\ell}2}
\int\Big(\int\Big|\int|\cT_m^{x',y'}f^{y'}(x'')| dy'\Big|^2 dx''\Big)^{1/2} dx'
\\
&\lc\delta_Q^{\ell/2}
\int\int\Big(\int
|\cT_m^{x',y'}f_Q^{y'}(x'')|^2 dx''\Big)^{1/2} dx'dy'
\\
&\lc\delta_Q^{\ell/2}
\int\int
\frac{2^{m(d-\ell-\s)}}{(1+2^{m}|x'-y'|)^N}\Big(\int_{y''}
|f_Q^{y'}(y'')|^2 dy''\Big)^{1/2} dx'dy'
\\
&\lc 2^{-m\s}\delta_Q^{\ell/2} 
\int_{x'}
\Big(\int_{y''}
|f_Q^{x'}(y'')|^2 dy''\Big)^{1/2} dx'
\\
&\lc 2^{-m\s} \delta_Q^{d/2} \|f_Q\|_2\lc 2^{-m\s}.
\tag 6.6
\endalign
$$

On the complement of $W_Q$  we use the kernel estimates of Lemma 4.3.

We split $a_m=\sum_{k\ge m-1} a_{k,m}$ where the kernel $K_{k,m}$ of 
$\cT[a_{k,m}]$ satisfies the estimate (4.7) with $\rho=0$. Consequently
since $|x''-\fS(y,x')|\approx |y''-S(x,y')|$ we have 
$$\int_{W_Q^c}|\cT_m^k  f_Q(x)| dx\lc 
4^{-k}\delta_Q^{-1} 2^{-m\s}\|f_Q\|_1\quad\text{ if }
 4^k\delta_Q\ge 1.
\tag 6.7$$

From  the 
gradient  estimates in  (4.8) and by using 
 the cancellation property of the atom $f_Q$ 
we get
$$\int|\cT_m^k  f_Q(x)| dx\lc 4^k\delta_Q 2^{-m\s}\|f_Q\|_1\quad\text{ if }
 4^k\delta_Q\le 1,\tag 6.8  $$ 
and the asserted  $H^1\to L^1$ bound follows from (6.6), (6.7) and (6.8).
\qed
\enddemo

\proclaim{Corollary 6.5}
Suppose that $0<\rho<(d-\ell)/2$ and $\sigma>2\rho$.  
Then  $T\in \cI^{\rho,-\sigma}$ is bounded on 
$L^{\frac{d-\ell}{d-\ell-\rho}}$ and bounded on $L^{\frac{d-\ell}{\rho}}$.
\endproclaim

\demo{\bf Proof} We shall prove  the
$L^{\frac{d-\ell}{d-\ell-\rho}}$ 
boundedness; by  \S2.3
this also implies the 
$L^{\frac{d-\ell}{\rho}}$ boundedness.

Let $a\in S^{\rho,-\s}$ and let $a_m$ be as in (6.5). Define
$$a_{m,z}(x,y,\tau,\xi)=
a_{m}(x,y,\tau,\xi)(1+|\tau|^2+|\xi|^2)^{(\rho_0(1-z)+\rho_1 z-\rho)/2}
(1+|\xi|^2)^{(\sigma-\sigma_0(1-z)-\sigma_1 z)/2}
$$
where $\sigma_0=\frac{d-\ell}{d-\ell-2\rho}(\sigma-2\rho)$, $\sigma_1=d-\ell$,
$\rho_0=0$ and $\rho_1=(d-\ell)/2$. Then $a_{m,\theta}=a_m$ for 
$\theta=2\rho/(d-\ell)$.
For $ \Re(z)=0$ the symbol
$a_{m,z}$ belongs to $S^{\rho_0,\sigma_0}$
and 
for $ \Re(z)=1$ it
 belongs to $S^{\rho_1,\sigma_1}$.
By Lemma 6.3 and Lemma 6.4 the operator
$\cT[a_{m,z}]$ is bounded from $H^1$ to $L^1$, with norm $(1+m)2^{-m\sigma_0}$
 if $\Re(z)=0$. By Lemma 6.1 and Lemma 6.2 it is bounded on $L^2$ with norm $O(1)$ if $\Re(z)=1$.
By interpolation we find that $\cT[a_{m}]$ is bounded on 
$L^{\frac{d-\ell}{d-\ell-\rho}}$ with norm
$O((1+m)2^{-m\s_0(1-\theta)})=O((1+m)2^{-m(\s-2\rho)})$. The assertion follows by summing in $m$.\qed
\enddemo

\head{\bf 7. $\boldkey L^{\boldkey p}\boldsymbol \to \boldkey L^{\boldkey q}$ estimates for Fourier integral operators}
\endhead

We begin by giving a different formulation of parts 1.2.1 and 1.2.2 of Theorem 1.2.
Suppose that  $0<\rho<(d-\ell)/2$ and $2\rho<\sigma<d-\ell$. Then statement 
1.2.1 of Theorem 1.2 says that
$T\in \cI^{\rho,-\sigma}$ maps $L^p\to L^q$ if $(1/p, 1/q)$ belongs to the closed trapezoid
with corners 
$(\frac{\r}{d-\ell}, \frac{\r}{d-\ell})$,
$(\frac{d-\ell-\r}{d-\ell},\frac{d-\ell-\r}{d-\ell})$,
$(1/p_{\rho,\sigma}, 1/q_{\rho,\sigma}),
(1/{q}_{\rho,\sigma}', 1/{p}_{\rho,\sigma}') $
where
$$
\aligned
&\frac{1}{p_{\rho,\sigma}}= \frac{d-\ell-\rho}{d-\ell}
-\frac{(\sigma-2\rho)\ell}{(d+\ell)(d-\ell)}
\\
&\frac{1}{q_{\rho,\sigma}}= \frac{d-\ell-\rho}{d-\ell}
-\frac{(\sigma-2\rho)d}{(d+\ell)(d-\ell)}.
\endaligned
\tag 7.1
$$

Observe that
$$
\frac{1}{p_{\rho_0,\sigma_0}}=
\frac{1}{q_{\rho_0,\sigma_0}}=\frac 12\quad\text{ if } \quad\rho_0=\frac{d-\ell}2,\quad 
\sigma_0=d-\ell,
\tag 7.2
$$
and if 
$$
\rho_1=0,\quad \sigma_1=(\sigma-2\rho)\frac{d-\ell}{d-\ell-2\rho}, \quad
\theta=\frac{d-\ell-2\rho}{d-\ell}
\tag 7.3
$$
then  $2\rho<\sigma<d-\ell$ implies $0<\sigma_1<d-\ell$ and we compute that 
$$
(1-\theta)\big(\frac{1}{p_{\rho_0,\sigma_0}},
\frac{1}{q_{\rho_0,\sigma_0}}\big)
+
\theta\big(\frac{1}{p_{\rho_1,\sigma_1}},
\frac{1}{q_{\rho_1,\sigma_1}}\big)=
\big(\frac{1}{p_{\rho,\sigma}},
\frac{1}{q_{\rho,\sigma}}\big).
\tag 7.4
$$
Therefore, one would like to prove Theorem 1.2 by 
interpolation from 
an  $L^{p_1}\to L^{q_1}$ 
result for operators in $\cI^{0,-\sigma_1}$ (already proved only for the case of 
 weakly singular Radon transforms)  
and an $L^2$ result for operators in
$\cI^{\frac{d-\ell}2, \ell-d}$. Unfortunately, 
operators in the 
latter class may fail to be bounded on $L^2$;  
this  somewhat complicates the interpolation  argument.

Performing a finite 
 finite conic partition of unity in the $\tau$ variables we may assume that 
$$\supp a\subset\big \{(x,y,\tau,\xi):
 |x|+|y|\le \eps^{10}, |\tau|+|\xi|\ge 2^{M+10} , \big|\tfrac{\tau}{|\tau|}-\vth|\le \eps+|\tau|^{-1}
\big\},
$$
for some given unit vector $\vth$ in $\bbR^\ell$, and $M$ is chosen as in \S5.

We shall now set up the various interpolation arguments.  We fix $\rho$ and $\sigma$
and use the abbreviation 
$$(p,q)=(p_{\rho,\sigma}, q_{\rho,\sigma} ),
\quad
(p_i,q_i)=(p_{\rho_i,\sigma_i}, q_{\rho_i,\sigma_i} ),
\quad i=1,2.
$$

We may split $T=T_{FIO}+T_{PsDO}$ where $T_{FIO}$ 
corresponds to a symbol which is supported where $|\tau |^{1/2}\ge|\xi|/2+2^{M+5}$ 
and
$T_{PsDO}$ corresponds to a symbol supported  in the complementary region.
Thus $T_{PsDO}=\cT[b]$ where
$b$ vanishes if  $|\tau |^{1/2}\ge 2|\xi|+10$ .
Let $$
W_z(\xi,\tau)=
(1+|\tau|^2+|\xi|^2)^{(\rho_0(1-z)+\rho_1 z-\rho)/2}
(1+|\xi|^2)^{(\sigma-\sigma_0(1-z)-\sigma_1 z)/2}
$$ and  $b_z(x,y,\tau,\xi)= b(x,y,\tau,\xi)W_z(\xi,\tau)$, so that 
$W_\theta=1$.  By Lemma 6.1 the operator  $\cT[b_z]$ is bounded on $L^2$
 if $\Re(z)=0$ and by Proposition 4.4  it is bounded from
$L^{p_1} $ to $ L^{q_1}$ if $\Re(z)=1$; all bounds are of admissible growth in $z$.
Thus $T_{PsDO}$ maps $L^p$ to $L^q$ by analytic interpolation.

Now we consider $T_{FIO}=\cT[a]$ where 
$a$ vanishes if  $|\tau |^{1/2}\le \max\{ 2^M, |\xi|/2\}$.
We first split off another  operator which behaves like $T_{PsDO}$.
Let $a_z=aW_z$ and $a_{k,m,z}=\beta_{k,m} a_z$ where $\beta_{k,m}$ is as in 
(4.5). Also let
$$\align
&a_{k,m,j,z} (x,y,\tau,\xi)
=a_{k,m,z} (x,y,\tau,\xi)\zeta(2^j|x'-y'|)
\\
&\widetilde a_{k,m,z} (x,y,\tau,\xi) 
=a_{k,m,z} (x,y,\tau,\xi) \zeta_0(2^k|x'-y'|)
\endalign
$$
Let
$$
\cV_{s,z}:= \sum_{k\ge s} \cT[\widetilde a_{k, k-s,z}]
$$

By Lemma 6.2 (i), with the choice $m(k)=k-s$,  the operator $\cV_{s,z}$ is bounded on 
$L^2$, uniformly in $s$,
if $\Re(z)=0$. By  Lemma 4.5 it is bounded from
$L^{p_1}\to L^{q_1}$ if $\Re z=1$; the bound is 
 $O(2^{-s(d-\ell-\sigma_1)})$;
 all bounds are admissible in $z$. 
Interpolating we see that $\cV_{s,\theta}$ maps $L^p\to L^q$ with norm
  $O(2^{-s(d-\ell-\sigma_1)\theta})
=O(2^{-s(d-\ell-\sigma)})$;  hence 
$\sum_{k,m} \cT[\widetilde a_{k,m}]$ maps $L^p$ to $L^q$.

It remains to estimate the operator 
$\sum_{k>0}\sum_{m<k}\sum_{j<k} \cT[a_{k,m,j,z}]$. We wish to use an angular 
Littlewood-Paley decomposition as in the proof of
Proposition 5.3.
Given a unit vector $v$ in $\bbR^{d-\ell}$ we make an angular localization in
$x'-y'$. By employing  a finite partition of unity it then suffices to bound 
$\sum_{k>0}\sum_{m<k}\sum_{j<k} \cT[\alpha_{k,m,j,z}]$
where
$$\alpha_{k,m,j,z}(x,y,\tau,\xi)=  
a_{k,m,j,z}(x,y,\tau,\xi)\zeta_0(\eps^{-5}|\frac{x'-y'}{|x'-y'|}-v|).$$
We choose $u$ as in (5.6) 
and  perform the change of variable $w\mapsto (w',w''+F(w;u))\equiv \cQ(w)$
 in \S 2.2, 
and define
$\fQ h(z)=h(\cQ (z))$.

As a result we have to show the $L^p\to L^q$ bound for the operator
$$\sum_{k>0}\sum_{m<k}\sum_{j<k} \fQ\cT[\alpha_{k,m,j,z}]\fQ^{-1}
=\sum_{k>0}\sum_{m<k}\sum_{j<k} \cT_{k,m,j}^z
\tag 7.5
$$
which has kernel
$$\sum_{k>0}
\sum_{m<k}\sum_{j<k} 
\iint e^{i[\inn{\tau }{y''-\widetilde S(x,y')}+\inn{x'-y'}{\xi}]} 
\widetilde \alpha_{k,m,j,z}(x,y,\tau,\xi)\, d\tau d\xi  
$$
where
$\inn {u}{\widetilde S_{x'}(x,x')}=
\inn {u}{\widetilde S_{y'}(x,x')}=0$ and 
$\widetilde \alpha_{k,m,j,z}(x,y,\tau,\xi)=
\alpha_{k,m,j,z}(\cQ(x),\cQ(y),\tau,\xi)\,g(x)/g(w)$,
and $g$ is smooth and positive.

We now use a Littlewood-Paley operators $L_k$ defined by
$L_k =\sum_{i=-4}^{4} \om(4^{-k+i}|D''|)$
and also the angular the Littlewood-Paley operator $P_{k,j}$ defined in (5.8).
Let
$$\cT_{k,m,j}=\cT[\widetilde \alpha_{k,m,j,\theta}].$$
We split 
$$
\sum_{k,m,j} \cT_{k,m,j}= 
\sum_{k,m,j} L_k\cT_{k,m,j}L_k
+\sum_{k,m,j} (I-L_k)\cT_{k,m,j}L_k
+\sum_{k,m,j} \cT_{k,m,j}(I-L_k)
$$
and then
$$
\sum_{k,m,j} L_k\cT_{k,m,j}L_k =(I+II)+(III+IV)+(V+VI)
$$
where
$$
\align
I+II&= \big[\sum\Sb k,m,j \\ m\le j\endSb+
\sum\Sb k,m,j \\ m> j\endSb\big]
L_k P_{k,j}\cT_{k,m,j}P_{k,j}L_k
\\
III+IV&= \big[\sum\Sb k,m,j \\ m\le j\endSb+
\sum\Sb k,m,j \\ m> j\endSb\big]
L_k (I-P_{k,j})\cT_{k,m,j}P_{k,j}L_k
\\
V+VI&=
 \big[\sum\Sb k,m,j \\ m\le j\endSb+
\sum\Sb k,m,j \\ m> j\endSb\big]
L_k \cT_{k,m,j}(I-P_{k,j})L_k.
\endalign
$$
We then split $I=\sum_{s\ge 0}I_s$ by linking $m=j-s$ for $s\ge 0$ and prove bounds 
for the expressions $I_s$ which decay in $s$. Similarly we split $II$ setting
$j=m-s$. The expressions $III, IV, V, VI$  are split into a double series depending on nonnegative parameters $r$, $s$; we prove then decay in $r,s$. We set $j=k-r, m=k-r-s$ when estimating 
$III$ and $V$ and  $j=k-r-s$, $m=k-r$ when estimating $IV$ and $VI$. In the
 following proposition  we state the relevant estimates for the pieces.

\proclaim{Proposition 7.1} Let $0\le \rho<(d-\ell)/2$ and $2\rho<\s<d-\ell$ and let
 $p=p_{\rho,\sigma}$, $q=q_{\rho,\sigma}$.
 There is $\delta=\delta(\rho,\sigma)>0$ so that the following estimates hold.

(i) For $s\ge 0$
$$
\Big\|\sum_{k>s}\sum_{s\le j<k} 
L_k P_{k,j}\cT_{k,j-s,j}P_{k,j} L_k\Big\|_{L^p\to L^q}\lc 
2^{-s\delta}
\tag 7.6
$$

(ii) For $s\ge 0$ 
$$
\Big\|\sum_{k>s}\sum_{s\le m<k}
L_k P_{k,m-s}\cT_{k,m,m-s}P_{k,m-s} L_k\Big\|_{L^p\to L^q}\lc 
2^{-s\delta}
\tag 7.7
$$

(iii)  For $s\ge 0$, $r\ge 0$,
$$
\align
&\Big\|\sum_{k>s+r}
L_k(I- P_{k,k-r})\cT_{k,k-r-s,k-r}P_{k,k-r} L_k\Big\|_{L^p\to L^q}\lc 
2^{-(r+s)\delta}
\tag 7.8
\\
&\Big\|\sum_{k>s}
L_k\cT_{k,k-r-s,k-r}(I-P_{k,k-r}) L_k\Big\|_{L^p\to L^q}\lc 
2^{-(r+s)\delta}.
\tag 7.9
\endalign
$$

(iv) For $s\ge 0$, $r\ge 0$,

$$\align
\Big\|\sum_{k>s+r}
L_k (I- P_{k,k-r-s})\cT_{k,k-r,k-r-s}P_{k,k-r-s} 
L_k\Big\|_{L^p\to L^q}\lc 2^{-r-s}
\tag 7.10
\\  \Big\|\sum_{k>0}
L_k \cT_{k,k-r,k-r-s}(I-P_{k,k-r-s}) L_k\Big\|_{L^p\to L^q}\lc
2^{-r-s}.
\tag 7.11
\endalign
$$

(v) For $j<k$, $m<k$
$$
\align &\|(I-L_k)\cT_{k,m,j}L_k\|_{L^p\to L^q}\lc 2^{-k}  
\tag 7.12
\\
&\|\cT_{k,m,j}(I-L_k)\|_{L^p\to L^q}\lc 2^{-k}  .
\tag 7.13
\endalign
$$
\endproclaim

Taking Proposition 7.1 for granted  we can complete the 
\demo{\bf Proof of Theorem 1.2}
Let
$p_{\r,\s}$ and $q_{\r,\s}$ be as in (7.1).
A combination of the estimates in Proposition  7.1 shows that the operator in (7.5) is bounded
 from
$L^{p_{\rho,\sigma}}$ to $L^{q_{\r,\s}}$. Together with the discussion preceding (7.5)
this yields
the $L^{p_{\rho,\sigma}}\to L^{q_{\r,\s}}$ bound  of the operator $\cT[a]$ 
where $a\in S^{\rho,-\s}$.
If we apply  this to  the adjoint operator we obtain  
the $L^{q_{\rho,\sigma}'}\to L^{p_{\r,\s}'}$ bound.
If $\rho>0$ we interpolate  with the $L^p\to L^p$ estimate in  \S6, and if $\rho=0$
we interpolate instead with the $H^1\to L^1$ bound in \S6. This yields
 the proof of statements 1.2.1 and 1.2.2. Statements 1.2.4 and  1.2.3 have already been proved in \S4 and \S5,
 respectively.\qed
\enddemo

We now give  a sketch of the

\subheading{Proof of Proposition  7.1}

We begin by estimating the main terms (7.6), (7.7) and use

\proclaim{Lemma 7.2}
Let $\cR^{\s_1}$ be as in (1.12) and  let $\Re(z)=1$.
Then
$$
|\cT^z_{k,m,j} f(x)|\lc \min\{ 2^{-(j-m)(d-\ell-\sigma_1)}, 2^{-(m-j)}\} 
\fM(\fQ\cR^{\sigma_1}\fQ^{-1}[f])
$$
where 
$\fM$ denotes the strong maximal function.
\endproclaim

\demo{Proof} 
This follows  from the kernel estimates (4.7) in a straightforward way. \qed
\enddemo

\demo{Proof of (7.6), (7.7)} By Theorem 5.1  we know that
$\cR^{\sigma_1}$ maps $L^{p_1}$ to $L^{q_1}$ and so does
$\widetilde{\cR}^{\sigma_1}\fQ\cR^{\sigma_1}\fQ^{-1}$.
Arguing as in the proof of Lemma 5.5, by the  Fefferman-Stein and 
Marcinkiewicz-Zygmund theorems
we therefore  have the vector-valued inequality
$$
\Big\|\Big (\sum_{j,k}|\fM \widetilde{\cR}^{\sigma_1}\!
 f_{j,k}|^2\Big)^{1/2}\Big\|_{q_1}
\lc
\Big\|\Big (\sum_{j,k}|f_{j,k}|^2\Big)^{1/2}\Big\|_{p_1}.
$$
We apply the $L^{q_1}\to L^{q_1}$ and 
$L^{p_1}\to L^{p_1}$ 
 Littlewood-Paley  inequalities for the Littlewood-Paley
 decompositions
$\{L_k P_{k,j}\}_{j,k}$ and Lemma 7.2 and obtain
$$
\Big\|\sum_{k>s}\sum_{s\le j<k} 
L_k P_{k,j}\cT^z_{k,j-s,j}P_{k,j} L_k\Big\|_{L^{p_1}\to L^{q_1}}\lc 
2^{-s(d-\ell)\sigma_1} 
\quad \text{ if }\Re(z)=1.
\tag 7.14
$$
By Lemma 6.2 and the almost  orthogonality of the Littlewood-Paley operators
$$
\Big\|\sum_{k>s}\sum_{s\le j<k} 
L_k P_{k,j}\cT^z_{k,j-s,j}P_{k,j} L_k\Big\|_{L^{2}\to L^2}\lc 1
\quad \text{ if }\Re(z)=0.
\tag  7.15
$$
(7.14) and (7.15)  prove (7.6) by interpolation and (7.7) is proved in the same way.
\enddemo

\demo{Proof of (7.8), (7.9), (7.10), (7.11)}
We analyze the kernel 
 of
$L_k (I-P_{k,j}) T^z_{k,m,j}$ which is given by
$$
\iint \iiiint  
e^{\ic \psi(x,t,h'',y,,\lambda,\eta'',\tau,\xi) }
\gamma_{k,m,j,z}(x,t,h'',y,\lambda,\eta,\tau,\xi)  
\,d\eta'' d\la d\tau d\xi  \,
dt dh'' 
$$
where
$$
\psi(x,t,h'',y,\lambda,\eta'',\tau,\xi) =-t\la-\inn{h''}{w''}+\inn{\tau }
{y''-\widetilde S(x'+tu,x''+h'',y')}+\inn{x'+tu-y'}{\xi}
$$
and
$$
\multline
\gamma_{k,m,j,z}(x,t,h'',y,,\lambda,\eta'',\tau,\xi) 
\\=
\sum_{i_1=-4}^4\zeta(4^{-k+i}|\eta''|)
\sum
\big(1-\sum_{i_2=-M}^M
\zeta(2^{-2k+j+i_2} |\lambda|) \big)
\widetilde \alpha_{k,m,j,z}(x'+tu,x''+h'',y,\tau,\xi).
\endmultline
$$
Arguing as in \S5 we first  integrate by parts with respect to $t$.
This yields the pointwise estimate 
$$\multline
2^{(2j-2k)N_2}
\int_w \frac{2^{2k-j}}{(1+2^{2k-j}|t|)^N}
\frac{2^{2k\ell}}{(1+2^{2k}|h''|)^N}
\\
\times \chi_j(x'+tu-y')
\frac{2^{m(d-\ell-\s)}}{(1+2^m|x'+tu-y'|)^N}
\frac{2^{2k\ell}}{(1+2^{2k}|y''-\widetilde S(x'+tu, x''+h'',y')|)^{N_1}} 
dt dh''
\endmultline
$$
here $N_2\gg N_1, N$ and
$\chi_j$ is the characteristic function of 
$\cup_\pm \pm[2^{-j-1},2^{-j+1}]$. A somewhat lengthy  but straightforward
calculation similar to the one 
for the term $\widetilde \cE^{1,i}$   in \S5 shows that for $s\le j\le k$
$$|L_k(I- P_{k,j})\cT^z_{k,j-s,j}f(x)|
\lc\int
 4^{j-k} 2^{-s(d-\ell-\s_1)}(|x'-y'|+|y''-\widetilde S(x,y')|^{1/2})^{\sigma_1-d-\ell}
|f(y)| dy, \quad \Re(z)=1,
$$
if $|x|\le \eps$ and better (trivial) decay 
estimates for $|x|\ge \eps$.

By using 
the $L^{p_1}\to L^{q_1}$ mapping property of the standard fractional integral operator and its vector-valued extension, together with the
$L^p$ inequalities for the Littlewood-Paley operator defined by $L_k$
(or $\widetilde L_k$ with  $\widetilde L_k L_k=L_k$) we obtain the estimate
$$\Big\|\sum_{k>s+r}
L_k(I- P_{k,k-r})\cT^z_{k,k-r-s,k-r}P_{k,k-r} L_k\Big\|_{L^{p_1}\to L^{q_1}}
\lc 
2^{-r} 2^{-s(d-\ell-\sigma)},\quad\Re(z)=1.
$$
By Lemma 6.2,
$\cT^z_{k,k-r-s,k-r}$ is bounded on $L^2$ if $\Re(z)=0$, uniformly in $s$, $r$ and $k$, and
 by the almost orthogonality of the $L_k$ (or $\widetilde L_k$) we get
$$\Big\|\sum_{k>s+r}
L_k(I- P_{k,k-r})\cT^z_{k,k-r-s,k-r}P_{k,k-r} L_k\Big\|_{L^{2}\to L^{2}}
\lc 1, \quad
\Re(z)=0.
$$ 
Analytic interpolation yields (7.8). The estimates (7.9), (7.10) and (7.11)
are proved in the same way.
\enddemo

\demo{Proof of (7.12), (7.13)} 
One writes out the integrals defining the kernels of the decompositions of
$L_lT^z_{k,m,j}$ and, if $|l-k|>2$ one gains factors $\min\{2^{-kN},2^{-lN}\}$
by integrating in the ${}''$-variables.\qed
\enddemo


\Refs
\ref\no 1 \by P. Brenner \paper $L_p$-$L_{p'}$ estimates for Fourier
integral operators related to hyperbolic equations \jour Math. Z.
\vol 152 \yr 1977 \pages 273--286 \endref




\ref\no 2 \by H. Carlsson, M. Christ, A. C\'ordoba, J. Duoandikoetxea,
J.L. Rubio de Francia, J. Vance, S. Wainger and D. Weinberg
\paper $L^p$ estimates for maximal functions 
and Hilbert transforms along flat curves in
$\Bbb R^2$
\jour Bull. Amer. Math. Soc. \vol 14\yr 1986\pages 263--267
\endref

\ref \no 3\by M. Christ
\paper Endpoint bounds for singular fractional integral operators
\jour preprint 1988
\endref

\ref\no 4\bysame\paper
Failure of an endpoint estimate for integrals along curves
\inbook Fourier analysis and partial differential equations
\bookinfo  ed. by J. Garcia-Cuerva, E. Hernandez, F. Soria and J. L. Torrea
\publ CRC Press \yr 1995
\endref


\ref\no 5 \by S. Cuccagna \paper Sobolev estimates for fractional and singular Radon transforms
\jour J. Funct. Anal. \vol 139\yr 1996\pages 94--118\endref





\ref\no 6\by G. Folland\book Real Analysis, modern techniques and their applications\publ Wiley\yr 1984\endref

\ref\no 7\by J. Garcia-Cuerva and J.-L. Rubio de Francia\book Weighted norm
inequalities and related topics
\publ North-Holland\yr 1985
\endref

\ref \no 8 \by I.M. Gelfand and G.E. Shilov \book Generalized functions Vol. I
\publ Academic Press \publaddr New York \yr 1964 \endref


\ref\no 9\by L. Grafakos \paper Strong type bounds for analytic families of
fractional integrals\jour Proc. Amer. Math. Soc.\vol 117\yr 1993\pages 653--663\endref

\ref\no 10 \by   A. Greenleaf, A. Seeger  and S. Wainger\paper
On X-ray transforms for 
rigid line complexes and 
integrals over curves in
$\Bbb R^4$\jour Proc. Amer. Math. Soc.\vol 127\yr 1999\pages 3533-3545\endref

\ref \no 11\by A. Greenleaf and G. Uhlmann
\paper Estimates for singular Radon transforms and pseudo-differential
operators with singular symbols
\jour J. Funct. Anal. \vol 89 \yr 1990\pages 202--232
\endref

\ref\no   12  \by L. H\"ormander \paper Fourier integral operators I \jour Acta Math.
\vol 127 \yr 1971 \pages 79--183 \endref


\ref\no 13\by R. Melrose \paper Marked Lagrangian Distributions
\jour manuscript \yr 1990\endref

\ref\no  14 \by A. Nagel, E. M. Stein and S. Wainger
\paper Differentiation in lacunary directions
\jour Proc. Nat. Acad. Sc. USA \vol 75\yr 1978\pages 1060--1062
\endref


\ref\no 15\bysame
\paper Hilbert transforms and maximal functions
related to variable curves
\inbook Harmonic Analysis in Euclidean spaces
(Proc. Sympos. Pure Math. Williams Coll., Williamstown,
Mass., 1978), Part 2, \pages 175--177\bookinfo
Proc. Sympos. Pure Math. XXXV\publ Amer. Math. Soc.\yr 1979
\endref




\ref\no 16 \by D. H. Phong and E.M. Stein
\paper Hilbert integrals, singular integrals and Radon transforms, I
\jour Acta Math. \vol 157\pages 99--157\yr 1986
\endref


\ref\no 17\bysame
\paper Singular Radon transforms and oscillatory integrals
\jour Duke Math. J.\vol 58\yr 1989\pages 347--369\endref

\ref\no 18\by F. Ricci and E. M. Stein
\paper Harmonic analysis on nilpotent groups and singular
integrals III: Fractional integration along manifolds
\jour J. Funct. Anal. \vol 86\yr 1989\pages 360--389
\endref



\ref \no 19 \by A. Seeger
\paper $L^2$ estimates for a class of singular oscillatory integrals
\jour Math. Res. Lett. \vol 1\yr 1994\pages 65--73
\endref


\ref\no 20 \by  A. Seeger and T. Tao\paper Sharp Lorentz space  estimates
for rough  operators\jour preprint\endref

\ref\no 21\by A. Seeger and S. Wainger
\paper
 Maximal and singular Radon transforms under convexity assumptions
 \jour in preparation\endref


\ref\no 22\by E.M. Stein\book Harmonic analysis: Real variable methods,
orthogonality and
 oscillatory integrals\publ Princeton Univ. Press \yr 1993
\endref




\endRefs
\enddocument